\documentclass[11pt]{article}
\usepackage{amssymb,amstext,amsmath,amsthm}

\newcommand{\nek}  {\newcommand}

\nek{\vyk}[1]{}

\nek{\ubf} {\fontseries{b}\selectfont}
\nek{\bfsl}{\bfseries\slshape}

\DeclareMathAlphabet{\skr}{U}{eus}{m}{n}{\skewchar\font'60}

\nek{\seci}  {\subsection}

\nek{\itsep}{\itemsep=0.4ex plus 0.15ex minus 0.1ex}
\nek{\tenu}[1]{
\def\theenumi{#1}\itsep}
\nek{\squarE}{\square}
\theoremstyle{plain}
\newtheorem{theorem}             {Theorem}[subsection]
\newtheorem{corollary}  [theorem]{Corollary}
\newtheorem{lemma}      [theorem]{Lemma}
\newtheorem{propo}      [theorem]{Proposition}

\newtheorem{thE}             {Theorem}

\theoremstyle{definition}
\newtheorem{problem}    [theorem]{Problem}
\newtheorem{numdef}     [theorem]{Definition}
\newtheorem{rem}        [theorem]{Remark}
\newtheorem{bag} [theorem]{Blanket agreement}
\newtheorem*{prF}{Proof}     
\nek{\bbag}{\begin{bag}}
\nek{\ebag}{\end{bag}}
\nek{\brem}{\begin{rem}}
\nek{\erem}{\qeD\end{rem}}
\nek{\bcor}{\begin{corollary}}
\nek{\ecor}{\end{corollary}}
\nek{\bdf} {\begin{numdef}}
\nek{\edf} {\qeD\end{numdef}}
\nek{\eDf} {\end{numdef}}
\nek{\ble} {\begin{lemma}}
\nek{\ele} {\end{lemma}}
\nek{\bpr} {\begin{propo}}
\nek{\epr} {\end{propo}}
\nek{\bqu} {\begin{problem}}
\nek{\equ} {\end{problem}}
\nek{\bte} {\begin{thE}}
\nek{\ete} {\end{thE}}
\nek{\bpF} {\begin{prF}}
\nek{\epF} {\qeD\end{prF}}
\nek{\epG} [1]{\qeG{#1}\end{prF}}
\nek{\qeD} {\hfill{$\qed$}}
\nek{\qeG} [1]
{\hfill{\hbox{$\qed$~({\sl#1\hspace{0.2ex}})}}}
\nek{\envur}[1]
{\begin{equation}\itsep#1\end{equation}} 
\nek{\ben}{\begin{enumerate}\itsep} 
\nek{\een}{\end{enumerate}} 
\nek{\bit}{\begin{itemize}\itsep}
\nek{\eit}{\end{itemize}}
\nek{\bay}{\begin{array}}
\nek{\eay}{\end{array}}
\nek{\slsP}{\hspace*{0.3ex}}
\nek{\rmsP}{\hspace*{0.5ex}}

\nek{\etc}{\hbox{{\sl etc}}}
\nek{\ie} {\hbox{{\sl i.}\slsP {\sl e.}}}
\nek{\ea} {\hbox{{\sl et.}\slsP {\sl al.}}}
\nek{\eg} {\hbox{{\sl e.}\slsP {\sl g.}}}
\nek{\wrt}{\hbox{w.\rmsP r.\rmsP t.}}
\nek{\noo}{\hbox{w.\rmsP l.\rmsP o.\rmsP g.}}
\nek{\Sat} {{\sf Saturation}}
\nek{\ttbox}[1] {{\mathtt{#1}}}
\nek{\dom}  {\mathop{\ttbox{dom}}}
\nek{\ran}  {\mathop{\ttbox{ran}}}
\nek{\st}   {\mathop{\ttbox{st}}}
\nek{\lh}   {\mathop{\ttbox{lh}}}
\nek{\al}{\alpha}
\nek{\da}{\delta}
\nek{\Da}{\Delta}
\nek{\ba}{\beta}
\nek{\vt} {\vartheta}
\nek{\ga}{\gamma}
\nek{\kpa}{\kappa}
\nek{\om} {\omega}
\nek{\Om} {\Omega}

\nek{\bbb}{\hspace{0.1pt}}
\nek{\dvoj}[1]{{\bbb{\mathbb #1}\bbb}}
\nek{\dN}{{\dvoj N}}
\nek{\dR}{{\dvoj R}}
\nek{\dS}{{\dvoj S}}
\nek{\adR} {{\mathord{\upa\dR}}}
\nek{\adN} {{\mathord{\upa\dN}}}
\nek{\adS} {{\mathord{\upa{\hspace{0.1ex}}\dS}}}
\nek{\gotsp}{\hspace*{0.5pt}}
\nek{\got}[1]{\mathord{\gotsp\mathfrak #1\gotsp}}
\nek{\cont}{\text{\mtho\large${\got c}$}}
\nek{\sta}{\mathord{\kern 0.1ex\vphantom{|}^\ast\kern -0.3ex}}

\nek{\sups}[2]
{\mathord{\kern 0.05em\vphantom{X}^{#2}\kern -0.17em #1}}
\nek{\upa}[1] {{\sups{#1}\ast}}

\nek{\TS} {\textstyle} 
\nek{\sq}  {\subseteq}
\nek{\kaz} {\forall\,}
\nek{\sus} {\exists\,}
\nek{\res} {{\hspace{0.05ex}\restriction\hspace{0.05ex}}}
\nek{\dm}  {$$}
\nek{\eqv} {\Longleftrightarrow}
\nek{\ti}  {\times}
\nek{\nin} {\not\in}
\nek{\imp} {\Longrightarrow}
\nek{\mpi} {\Longleftarrow}
\nek{\bez} {\setminus}
\nek{\lra} {\longrightarrow}
\nek{\iy}  {\infty}
\nek{\limp}{\,\imp\,}

\nek{\mto} {\mapsto}
\nek{\lmto}{\longmapsto}
\nek{\onto}{\stackrel{\rm onto}\lra}
\nek{\ang} [1]{\langle #1\rangle}
\nek{\ans} [1]{\{\hspace{0.02ex}#1\hspace{0.02ex}\}}
\nek{\sis}[2] {{\ans{#1}}\vphantom{{}^x_x}_{#2}} 
\nek{\dd} [1] {$\rsur\mtho#1\qsur$-}
\nek{\itla}{\item\label}

\nek{\skl}{\hbox{\mtho\large$($}} 
\nek{\skp}{\hbox{\mtho\large$)$}} 
\nek{\mtho}{\mathsurround=0mm}
 
\nek{\mthf}{\mathsurround=0.4ex}
\nek{\msur}{\hspace*{-1\mathsurround}}
\nek{\qsur}{\hspace{0.2\mathsurround}}
\nek{\rsur}{\hspace{0.4\mathsurround}}
\nek{\noi}{\noindent}
\nek{\vom}{\vspace{1mm plus 0.3mm minus 0.3mm}}
\nek{\vtm}{\vspace{2mm plus 0.6mm minus 0.6mm}}

\nek{\poq}{\underline}

\nek{\rE} {\qE} 
\nek{\rPi}{\mathrel{\sf\Pi}}

\nek{\aeq}{\simeq}
\nek{\obr}{^{-1}}

\newlength{\dxii}
\nek{\zo} {,\linebreak[0]}
\nek{\zi} {,\linebreak[0]\,}
\nek{\zd} {,\linebreak[0]\:}
\nek{\zt} {,\linebreak[0]\;}

\nek{\mul}[1] {\begin{multline}#1\end{multline}}%
\nek{\muL}[1] {\begin{multline*}#1\end{multline*}}%

\nek{\otn} [3] {#1\mathrel{#2}#3}
\nek{\dtn} [3] {{#2}(#1,#3)}

\nek{\fs}[2]{{\bf\Sigma}^{#1}_{#2}}
\nek{\fso} [1] {\fs0{#1}}

\nek{\hh} {{\skr H}}
\nek{\cA} {{\skr A}}
\nek{\cF} {{\skr F}}
\nek{\cN} {{\skr N}}
\nek{\ps} {{\skr P}}
\nek{\pv} [1] {\ps_{\mathtt{int}}(#1)}

\nek{\qE} {\mathbin{\sf E}}
\nek{\qV} {\mathbin{\sf VIT}}
\nek{\FD} {\mathbin{\sf F\hspace{-0.1ex}D}}
\nek{\Fd} [1] {\mathbin{\FD\res{#1}}}
\nek{\qF} {\mathbin{\sf F}}
\nek{\qR} {\mathbin{\sf R}}
\nek{\EO} [1] {\mathbin{\hspace{0.05ex}{\En}\res{#1}\hspace{0.05ex}}}
\nek{\qei}  [1] {\mathbin{{\sf M}_{#1}}}
\nek{\qej}  [1] {\mathbin{{\sf M}^{#1}}}
\nek{\qeip} [1] {\mathbin{{\sf M}'_{#1}}}
\nek{\qeu}  {\qei U}
\nek{\qeup} {\mathbin{\sf R}}
\nek{\qen} {\qei{\mathbb N}}

\nek{\rei} [1] {\mathbin{{\sf R}_{#1}}}
\nek{\reil}[1] {\rei{\log#1}}
\nek{\ru}  [1] {\mathbin{{\sf R}_{\ans{#1}}}}

\nek{\bine}[1] {\sg_{#1}}

\nek{\qd} {countably determined}
\nek{\Qd} {Countably determined}
\nek{\kd}[1] {{\bf CD}\hspace{0.2ex}[#1]}

\nek{\ddn} {\fuk 2\adN}%
\nek{\dn}  {\fuk 2\dN}%
\nek{\nn}  {\fuk\dN\dN}%
\nek{\dln} {\fuk 2{<\om}}
\nek{\nln} {\fuk\dN{<\om}}

\nek{\aci} {\addtocounter{enumi}}

\nek{\er} {ER}

\nek{\qdr} {\le_{\tt C\hspace{0.2ex}D}}
\nek{\qde}  {\equiv_{\tt C\hspace{0.2ex}D}}
\nek{\qdl}  {<_{\tt C\hspace{0.2ex}D}}

\nek{\bos} {B}
\nek{\bor} {\le_{\tt\bos}}
\nek{\boe}  {\equiv_{\tt\bos}}
\nek{\bol}  {<_{\tt\bos}}

\nek{\bore}  {\equiv_{\tt B}}

\nek{\cu} [1] {{\text{\large\boldmath\mtho$\sqcup$}\ans{#1}}} 
\nek{\ceu}[1] {\qei{\ans{#1}}}

\nek{\cv} [1] {{\text{\large\boldmath\mtho$\sqcap$}{\ans{#1}}}}
\nek{\cev}[1] {\qej{\ans{#1}}}

\nek{\qev} {\qei V}

\nek{\dens} {\mathop{\tt dens}}
\nek{\card} {\mathop{\tt card}}
\nek{\ima}  {\mathord{\hspace{0.4ex}\text{''}}}

\nek{\adens}{\mathop{\sta\text{\tt dens}}}

\nek{\rav} [1] {\mathop{\text{\sf D}(#1)}}
\nek{\Exp} [1] 
{\mathop{\text{\sf D}_{\text{\tt ext}}(\fuk2{#1})}} 

\nek{\api} {w}
\nek{\apI} {W}

\nek{\eqr} {equivalence relation}

\nek{\gU} 
{\text{\mtho\boldmath$\hspace{0.1ex}\mathfrak U\hspace{0.1ex}$}}

\nek{\fuk}[2] {\mathord{^{#2}\hspace*{-0.10ex}{#1}}}
\nek{\fuko}[2] {\mathord{^{#2}\hspace*{-0.00ex}{#1}}}
\nek{\fuki}[2]{\mathord{^{#2}\hspace*{-0.4ex}{#1}}}

\nek{\sui} [1] {{#1}{}_{\text{\tt int}}}
\nek{\sbi} [1] {\sui{(#1)}}
\nek{\ifuk}[2] {\sui{(\fuki{#1}{#2})}}

\nek{\dt} [1] {\fuk2{#1}}

\nek{\fse} {\adS}

\nek{\ups} [1] {y_{#1}}

\nek{\bok} [1] {\#_{\text{\tt B}}(#1)}

\nek{\dst} {descriptive set theory}

\nek{\ek} [1] {[#1]_{\qE}}
\nek{\pro} {\mathop{\ttbox{prfl}}}

\nek{\vpi} {\varphi}
\nek{\vep} {\varepsilon}
\nek{\Sg}  {\Sigma}
\nek{\La}  {\Lambda}

\nek{\leb} {\mathop{\underline{\text{\tt mes}}}}
\nek{\leB} {\mathop{\overline{\text{\tt mes}}}}

\nek{\dda} {\dd\ast}

\nek{\sg}{\sigma}
\nek{\sd}{\mathbin{\Delta}}

\nek{\pu}  {\emptyset}
\nek{\sneq}{\subsetneqq}
\nek{\dop} [1] {{\complement}{#1}}
\nek{\we}  {\mathbin{^\wedge}}

\nek{\cX} {{\skr X}}

\nek{\ZFC} {\text{\bf ZFC}}
\nek{\ZF}  {\text{\bf ZF}}

\nek{\dde}{\dd\qE}

\nek{\cut}{\mathop{\text{\tt cut}}}

\nek{\dZ}{{\dvoj Z}}

\nek{\cdc} {\mathop{\ddag}}
\nek{\Eo} {\qE_0}
\nek{\En} {\mathop{{\displaystyle\qE_0^{\text{\rm nst}}}}}
\nek{\Eoq}{\text{``\mtho$\Eo$''}}

\nek{\bad} {_{\text{\tt bad}}}
\nek{\kon} {_{\text{\tt fin}}}
\nek{\bes} {_{\text{\tt inf}}}

\nek{\lek} {\preccurlyeq}

\nek{\frat} [2] {{#1}-{#2}}
\nek{\frakk}[2] {{#1}-{#2}}

\nek{\doo} [1] {\Omega_{\ans{#1}}}

\nek{\tx} {\tilde x}
\nek{\tg} {\tilde\sg}

\nek{\resb} [2] {{#1}\res_{\ge#2}}

\nek{\geo} {} 

\nek{\cud} [1] {\cu{2^{#1}}}
\nek{\cvd} [1] {\cv{2^{#1}}}

\nek{\pol} {``Polish''}

\nek{\wU}{\widetilde U}

\nek{\ent}[1]{{\cal E}(#1)}

\title{Borel and countably determined  
reducibility in nonstandard domain}

\author{%
Vladimir Kanovei,\thanks{Contact author.}
\thanks{
Support of DFG acknowledged. 
Moscow Center for continuous mathematical education,
Bol.~Vlasevski 11, Moscow, 121002, Russia,
{\tt kanovei@math.uni-wuppertal.de}.}%
\and
Michael Reeken,
\thanks{Department of Mathematics, University of Wuppertal,
Wuppertal, 42097, Germany,
{\tt reeken@math.uni-wuppertal.de}.}%
}

\date{December 2001}

\begin{document}\maketitle

\begin{abstract}
We consider, in a nonstandard domain, reducibility of 
equivalence relations in terms of the Borel reducibility 
$\bor$ and the countably determined 
(CD, for brevity) reducibility ${\qdr}.$
This reveals phenomena partially analogous to those discovered  
in modern ``standard'' descriptive set theory. 
The \dd\qdr structure of CD sets (partially) and 
the \dd\bor structure of Borel sets (completely) 
in $\adN$ is described. 
We prove that all ``countable'' 
(\ie, those with countable equivalence classes) 
CD \eqr s (\er s) are CD-smooth, but not all are B-smooth: 
the relation $x\qen y$ iff $|x-y|\in\dN$ is a counterexample. 
Similarly to the Silver dichotomy 
theorem in Polish spaces, any CD equivalence relation on 
$\adN$ either has at most continuum-many classes 
(and this can be witnessed, in some manner, by a \qd\ function) 
or there is an infinite internal set of pairwise inequivalent 
elements. 
Our study of {\it monadic\/} \eqr s, \ie, those of the form 
$x\qeu y$ iff $|x-y|\in U,$ where $U$ is an additive \qd\ cut 
(initial segment) demonstrates that these \er s split in 
two linearly \dd\bor(pre)ordered families, associated with 
countably cofinal and countably coinitial cuts, and the 
equivalence $u\FD v$ iff $u\sd v$ is finite, on the set of 
all hyperfinite subsets of $\adN,$ \dd\bor reduces all 
``countably cofinal'' \er s but does not \dd\qdr reduce 
any of ``countably coinitial'' \er s.
\end{abstract}

\mthf

Classical \dst\ (DST, for brevity) is mainly
concentrated on sets in Polish (complete separable) 
spaces, see Kechris~\cite{dst}.
It was discovered in 80s that ideas of classical DST can be
meaningfully developed in a very different setting of
nonstandard analysis, where Polish spaces are replaced by
internal hyperfinite sets.
This alternative version of \dst\ is called
``hyperfinite'', or ``nonstandard'' DST.
It allows to define Borel and projective hierarchies of 
subsets of a fixed infinite internal 
(for instance, hyperfinite) 
domain in quite the same manner as \pol, \ie, classical 
DST does, but beginning with internal sets 
at the initial level rather than open sets. 
Generally, the structures studied by the ``nonstandard'' 
DST appear to be similar, in some aspects, to those  
considered in the \pol\ \dst, 
but different in some other aspects. 
As for the proofs, they are mainly based on
very different and rather combinatorial ideas, and 
(countable) \Sat, of course, see Keisler~\ea~\cite{kkml}.
``Nonstandard'' DST also involves objects which hardly have
any direct analogy in the \pol\ setting, 
like \qd\ sets, leading to a remarkably interesting
mixture of \pol\ and nonstandard concepts and methods.

This note is written in attempt to find nonstandard
analogs of concepts which attract a lot of attention in
\pol\ DST nowadays: the structure of definable
(usually, Borel or analytic) equivalence relations in terms 
of Borel (sometimes more complicated) reducibility 
of associated quotient structures.
Our results will be related to \qd\ (or CD), in particular, 
Borel sets and equivalence relations on $\adN$ and hyperfinite 
domains, and the reducibility by \qd, in particular, by 
Borel maps.

It is an important difference with the \pol\ DST that  
while classically all uncountable Polish spaces are Borel 
isomorphic, hence indistinguishable \wrt\ topics in Borel 
reducibility, in ``nonstandard'' setting any two infinite 
hyperfinite sets $X,\,Y$ admit a Borel bijection iff 
$\frac{\#X}{\#Y}\simeq 1$ and admit a CD bijection iff 
iff $\frac{\#X}{\#Y}$ is neither infinitesimal 
nor infinitely large (see Proposition~\ref{cdb} below). 
This makes the structure of CD equivalence relations 
dependent not only on their intrinsic nature, \ie, the 
method of definition, but also on the size of the domain, 
which can be any internal infinite hyperfinite set of $\adN.$
(However see the last remark in Section~\ref P.)

This effect shows up already at the level of 
{\ubf B-smooth} \er s 
(those which admit a Borel enumeration of equivalence 
classes), 
which leads us to the study of Borel sets in terms of 
the relation $X\bor Y$  
meaning the existence of a Borel injection $\vt:X\to Y.$  
We prove (Theorem~\ref{cardB}) that  
any Borel subset of $\adN$ admits a Borel bijection onto 
a Borel cut (that is, initial segment) in $\adN,$ therefore, 
two Borel sets are comparable via the existence of a Borel 
injection, and generally there is a comprehensive 
classification of Borel subsets of $\adN$ modulo 
$\boe$ (that is, {\ubf Borel cardinalities}). 

A complete classification of \qd\ sets modulo $\qde$ 
is not known, yet we show 
(Theorem~\ref{tCD}) that, for any CD set $X\sq\adN,$ 
\poq{either} there is a unique {\it additive\/} CD cut 
$C\sq\adN$ 
(which can be equal to $\dN$ or $\adN$ itself) 
with $X\qde C,$ \poq{or} there is a hyperinteger 
$c\in\adN\bez\dN$ such that $c/\dN\qdl X\qdl c\dN.$ 
As a matter of fact 
we don't know whether the \poq{or} case really takes place. 
Anyway, we prove (Theorem~\ref{rer-c}) that any CD set 
$X\sq c\dN$ with $c/\dN\qdr X$ satisfies $X\qde M,$ where 
$M$ is a union of {\it monads\/} -- sets of the 
form $x+(c/\dN)\zd x\in c\dN,$ but whether such a 
set can satisfy $c/\dN\qdl X\qdl c\dN$ is not known. 

In \pol\ theory, some most elementary examples of 
non-smooth (in the sense of Borel enumerations, of course) 
ERs belong to the type of {\ubf countable\/} ones, \ie, 
with all equivalence classes at most countable. 
We prove (Theorem~\ref{sq}) that, on the contrary, in the 
``nonstandard'' DST any countable \qd\ ER $\qE$ 
admits a CD {\it transversal\/},  
\ie, a set which has exactly one
common element with each \dde class, hence, is CD-smooth 
(but not necessarily has a Borel transversal and is 
B-smooth, \ie, with a Borel 
enumeration of the equivalence classes).
This generalizes a recent theorem of Jin~\cite j   
that the (countable) equivalence relation $\qen$ 
defined on $\adN$ by $x\qen y$ iff $|x-y|\in\dN$ 
admits a \qd\ transversal and is CD-smooth. 
On the other hand, by a typical measure-theoretic 
argument, $\qen$ is not Borel-smooth and does not admit 
a Borel transversal; this is a transparent demonstration 
of differences between Borel and \qd\ structures. 

Theorem~\ref{Ts} belongs to the category of 
{\ubf dichotomy theorems}: 
in particular (the actual result is more general), 
it asserts that a CD equivalence relation has at most 
\dd\cont many equivalence classes or else admits an 
infinite internal set
of pairwise inequivalent elements. 
This has obvious similarities with the known theorems 
of \pol\ \dst, saying that a coanalitic (Silver), resp.,
analytic (Burgess) \er\ on a Polish space has
$\le\aleph_0,$ resp., $\le\aleph_1$ equivalence classes,
or admits an uncountable closed set of pairwise inequivalent 
elements. 
Generally speaking, the cardinality of continuum cannot 
be improved to any smaller value in Theorem~\ref{Ts}, yet 
in the case of \er s of class $\Sg^0_1$ it can be replaced 
by $\aleph_0,$ which improves upon Henson's \cite{he} 
theorem that any \qd\ set either is countable or contains 
an infinite internal subset.

An important class of \qd\ \er s which contains mostly 
non-CD-smooth relations, is the class of  
{\ubf monadic} \eqr s. 
Given an additive cut (initial segment) $U\sq\adN,$ we 
define $x\qeu y$ iff $|x-y|\in U,$ for all $x,\,y\in\adN.$ 
Since any additive CD cut 
(with trivial exceptions of $\pu$ and $\adN$) 
is either countably cofinal or countably coinitial 
(\ie, of the form, resp., 
$\bigcup_n[0,a_n)$ or $\bigcap_n[0,a_n),$ 
where $\sis{a_n}{n\in\dN}$ is strictly increasing, resp., 
decreasing sequence of hyperintegers), 
\qd\ monadic \er s split into two distinct families 
of {\it countably cofinal\/} and {\it countably coinitinal\/} 
monadic \er s. 
 
Our study of the reducibility phenomena among monadic \eqr s 
in Sections \ref{mer} --- \ref{tfin} 
(summarized in Theorem~\ref T) 
shows that \er s are mutually \dd\bor comparable within each 
of these two families, in such a way that the direction of 
$\bor$ between two monadic \er s $\qeu,\,\qev$ is determined 
by the relative 
rate of growth of countable cofinal sequences in $U,\,V$  
(or, in the countably coinitial case, of coinitial sequences in 
$\adN\bez U\zt\adN\bez V$), rather than by the relative size of 
the cuts $U,\,V,$ 
moreover, the \dd\qdr structure within either of the two 
families (but not between them) coincides with the 
\dd\bor structure.
It turns out that, in each of the two families, there is a 
subclass of \dd\bor minimal (and \dd\qdr minimal) 
\er s, namely, those generated 
by cuts of the form $c\dN$ or $c/\dN\zt c\in\adN$ 
(in, resp., countably cofinal, coinitial case).  
Further, among all monadic \er s only those of the form 
$\qei{c\dN}$ are CD-smooth 
(and all of them even admit a CD transversal, essentially by 
Jin~\cite j), 
but none of them is Borel-smooth. 
In addition, there is no relationship, in terms of 
$\bor$ or $\qdr,$ between 
countably cofinal and countably coinitial \er s except that 
we have $\qei{c\dN}\qdr\qev$ for any countably coinitial 
\eqr\ $\qev.$ 

Finally, we show in Section~\ref{eo} that monadic \er s 
induced by countably cofinal cuts admit a natural 
upper \dd\bor bound, namely, the \eqr\ of equality 
of hyperfinite subsets of $\adN$ modulo a finite set. 
We denote this \er\ by $\FD;$ it has some analogy with 
the equivalence relation of equality of infinite subsets 
of $\adN$ modulo a finite set, extensively studied in 
\pol\ \dst. 
We prove that $\qeu\bol\FD$ holds for any countably 
cofinal additive cut $U$ but fails for any countably 
coinitial additive $U.$ 
It is not clear whether $\FD$ is a \poq{minimal} upper bound 
for countably cofinal monadic \er s: this and some other 
open problems are considered in the final Section~\ref P.

\seci{Notation}
\label N


$\fuki XY$ is the set of all functions $f:Y\to X,$ 
while $x^y$ will denote only the arithmetical power 
operation in standard and nonstandard domains.


$\dln=\bigcup_{n\in\dN}\fuk2n$ 
is the set of all finite binary sequences. 

$s\we a$ is the extension of a finite sequence $s$ by 
a new rightmost term $a.$

$\lh s$ is the length of a finite sequence $s.$ 

$f\ima X=\ans{f(x):x\in X\cap\dom \vt},$ the 
\dd f{\it image\/} of a set $X.$ 

If $P$ is a set of pairs then   
$\otn xPy$ and $\dtn xPy$ mean that $\ang{x,y}\in P.$\vom


{\ubf Nonstandard setup.} 
Some degree of the reader's acquaintance with basic notions
of ``hyperfinite'' descriptive set theory is assumed; we give 
\cite{kkml} as the basic reference.
All ``nonstandard'' notions below, for instance $\adN,$
are related to a fixed countably saturated nonstandard
universe $\gU,$ whose elements will be referred to as
{\it nonstandard\/} (internal or external) sets.

In the remainder, we typically use letters like 
$i,\,j,\,k,\,m,\,n$ (with indices) for elements of $\dN,$ 
and letters like $a,\,b,\,c,\,h,\,x,\,y,\,z$
for elements of $\adN.$

$\pv X$ is the set of all internal subsets of a nonstandard 
set $X.$ 

If $X,\,Y$ are internal sets then $\ifuk XY$ is the set of 
all {\it internal\/} $f:Y\to X.$ 

Numbers $c\in\adN$ (standard or nonstandard) will be 
systematically identified with the sets $[0,c)=\ans{x:x<c}$ 
of all smaller numbers. 
We shall often use $\dt c,$ instead of the more pedantical 
$\sbi{\fuk2{[0,c)}}$ to denote the (internal) set of all 
internal functions $\xi:c=[0,c)\to2$.

$\#X\in\adN$ is the number
of elements of a hyperfinite set~$X.$

Let $r\simeq q$ mean that the difference $r-q$ is 
infinitesimal. 
For any bounded hyperrational $\al$   
(\ie, $\al<c$ for some $c\in\dN$) there is a unique 
standard real number $r,$ denoted by $\st \al,$ the
{\it standard part\/} of $\al,$ such that $\al\simeq r.$ 
If $\al$ is unbounded then put $\st\al=+\infty$.
\vom

{\ubf Borel and \qd\ sets.} 
Classes $\Sg^0_1\zt \Pi^0_1$ consist of countable 
unions, resp., intersections of internal sets. 
{\it Borel\/} sets form the least \dd\sg algebra which 
contains all internal sets; for instance, all sets in 
$\Sg^0_1\cup\Pi^0_1$ are Borel.

Following Henson \cite{he}, sets of the form
\envur{
\tag{\dag}
\textstyle
\label{cd1}
X=\bigcup_{b\in B}
\skl
\bigcap_{m\in B}X_{m}\,\cap\,\bigcap_{m\nin B}\overline X_{m}
\skp\,,\,
\text{ \ \parbox[t]{0.40\textwidth}
{\hfill where all sets $X_m$ are internal,\\[0.4ex] 
\hspace*{-8ex}\hfill
$B\sq\ps(\dN),$ and $\overline X_m=\bigcup_nX_n\bez X_m$,}}
}
are called {\it \qd\/}, in brief CD. 
(Any reasonable version of this concept for Polish spaces
yields the collection of all sets of the space.)
There are several slightly different ways to define this
class of sets, for instance,
\envur{
\tag{\ddag}
\textstyle
\label{cd}
X=\bigcup_{f\in F}\bigcap_{m\in\dN}X_{f\res m}\,,\,\
\text{ \ \parbox[t]{0.55\textwidth}
{where all $X_s\zt s\in\dln,$ are internal, \\[0.3ex]
$F\sq\dn,$
and $X_t\sq X_s$ whenever $s\subset t$.}}
}
(See, \eg, \cite j. 
To convert \eqref{cd} to \eqref{cd1}, let $B$ consist of
all sets $b\sq\dln$ containing a subset of the form 
$\ans{f\res m:m\in\dN}\zt f\in F,$ and apply any bijection 
$\dln$ onto $\dN.$
To convert \eqref{cd1} to \eqref{cd}, put
$X_s=\bigcap_{k<m}X'_k$ for any
$s=\ang{i_0,...,i_{m-1}}\in\dln,$ where $X'_k=X_k$ whenever
$i_k=1$ and $X'_k=\adN\bez X_k$ otherwise, then let $F\sq\dn$ 
be the set of all characteristic functions of sets in $B.$)

All Borel sets are \qd, but not conversely.

A map is Borel, \qd\ if it has a Borel, resp., 
CD graph.\vom

{\ubf Cuts.} 
Initial segments of $\adN$ (including $\pu\zo\dN\zo\adN$) 
are called {\it cuts\/}.
A cut $U$ is {\it additive\/} if $x+y\in U$ whenever $x,y\in U.$
Given a CD cut $U,$ the sets
\dm
\textstyle
U\dN=\bigcup_{n\in\dN,\,x\in U}[0,xn]
\quad\text{ and }\quad
U/\dN=\bigcap_{n\in\dN,\,x\in U}[0,\frac xn]
\dm
are additive CD cuts, $U/\dN\sq U\sq U\dN,$
$U/\dN$ is the largest additive cut included in $U$ while
$U\dN$ is the smallest additive cut including $U.$
In particular, let $c/\dN=[0,c)/\dN$ and $c\dN=[0,c)\,\dN$ for
any $c\in\adN$.

If $U$ is an additive cut then $\log U=\ans{h:2^h\in U}$ is 
also a cut (not necessarily additive) and 
$U=2^{\log U}=\bigcup_{h\in\log U}[0,2^h)$.

Internal cuts are $\pu\zo\adN,$ and those of the form 
$c=[0,c)\zt c\in\adN.$ 
Non-internal cuts can be obtained with the following 
general procedure. 
If $\sis{a_n}{n\in\dN}$ is a strictly increasing, 
resp., decreasing sequence in $\adN$ then we define a 
{\it countably cofinal\/} cut $\cu{a_n}=\bigcup_n[0,a_n),$ 
resp., {\it countably coinitial\/} cut 
$\cv{a_n}=\bigcap_n[0,a_n).$  
Both types consist of Borel sets of classes resp.\ $\Sg^0_1$ 
and  $\Pi^0_1.$ 

Cuts of the form $c+\dN=\ans{c+n:n\in\dN}$ and 
$c-\dN=\ans{c-n:n\in\dN}\;\,(c\nin\dN)$ are countably 
cofinal, resp., coinitial, but not additive 
(unless $c\in\dN$ in $c+\dN$). 

\ble
\label{cdcut}
Any CD cut\/ $\pu\ne U\sneq\adN$ it either countably cofinal 
or countably coinitial or contains a maximal element 
(and then is internal).
\ele
\bpF
Let $U=\bigcup_{f\in F}\bigcap_{m\in\dN}X_{f\res m},$
where $F$ and the sets $X_s$ are as in \eqref{cd}.
Put $\cut X=\bigcup_{x\in X}[0,x]$ for any set $X\sq\adN,$
the least cut which includes $X.$
By \Sat,
$U=\cut U=\bigcup_{f\in F}\bigcap_{m}U_{f\res m},$
where
$U_s=\cut X_s,$ hence, $U_s=[0,\mu_s],$ where
$\mu_s=\max X_s\in\adN$ for all $s\in\dln.$
If there is $f\in F$ with $U=\bigcap_{m}U_{f\res m}$
then the sequence $\sis{h_{f\res m}}{m\in\dN}$ witnesses that
$U$ is countably coinitial, or contains a maximal element
if the sequence is eventually constant.
Otherwise, by \Sat, for any $f\in F$ there is $m_f\in\dN$ such
that $h_{f\res m_f}\in U.$
Let $S=\ans{f\res m_f:f\in F};$ this is a countable set and
easily $U=\bigcap_{s\in S}[0,\mu_s],$ so that $U$ is either
countably cofinal or contains a maximal element.
\epF

\seci{Equivalence relations and reducibility: preliminaries} 
\label{prel}

Suppose that $\qE,\,\qF$ are \qd\ equivalence relations
(\er s, for brevity) on (also \qd) sets $X,\,Y.$
We write $\qE\qdr\qF,$ in words: 
{\it$\rE$ is CD-reducible to\/ $\qF,$\/}
iff there is a CD map (called: {\it reduction\/})
$\vt:X\to Y$~\footnote
{\ To apply $\qdr$ to non-CD relations, we should have 
used the existence of a CD map $\vt$ with $X\sq\dom\vt$ 
and $\vt\ima X\sq Y,$ but we'll not consider anything 
more complicated than CD below, in fact, mainly Borel 
\er s will be considered.}
such that we have ${x\qE x'}\eqv{\vt(x)\qF\vt(y)}$
for all $x,\,x'\in X.$~\footnote 
{\ It would be not less reasonable, but obviously longer, 
to write ${X/{\qE}}\qdr{Y/{\qF}}.$} 
We write $\qE\qde\qF$ if both $\qE\qdr\qF$ and $\qF\qdr\qE,$
and $\qE\qdl\qF$ iff $\qE\qdr\qF$ but not $\qF\qdr\qE.$ 
Changing ``\qd'' and ``CD'' to ``Borel'' in these definitions, 
we obtain the relations ${\bor}\zd{\boe}\zd{\bol}$ of 
{\it Borel\/} reducibility. 

Informal meaning of $\qE\qdr\qF$ and $\qE\bor\qF$ is that
{\it $\qF$ has at least as many equivalence classes as\/
$\qE,$ and this is witnessed by a CD, resp., Borel map\/}.\vom

{\ubf Equalities, smooth \er s, transversals.} 
For any set $A,$ the {\it equality relation\/} 
$\rav A$ ($\text{\sf D}$ from ``diagonal'') 
is defined on $A$ by $x\rav A y$ iff $x=y.$ 
These are the simplest of \er s; 
in many aspects $\rav A$ can be identified with $A.$

Similarly to the \pol\ \dst, say that an \er\ $\qE$ on a 
set $X$ is {\it CD-smooth\/}, resp., {\it B-smooth\/}, 
if $\qE \qdr \rav\adN,$ resp., $\qE \bor \rav\adN,$ 
\ie, there is a \qd, resp., Borel map
$\vt,$ with $X\sq\dom\vt$ and $\ran\vt\sq\adN$ such that 
$x\qE x'$ iff $\vt(x)=\vt(x'):$
this means that \dd\qE classes admit a
CD enumeration by hyperintegers. 

A {\it transversal\/} of an \eqr\ $\qE$ is a set which has 
exactly one common element with each \dd\qE equivalence 
class. 
Easily any Borel \er\ $\qE$ on a set $X\sq\adN,$ having a 
Borel transversal $W\sq X,$ is B-smooth: 
let $\vt(x)$ be the only element of $W$ equivalent to $x.$ 
Similarly any CD \eqr\ $\qE$ having a CD transversal is CD-smooth.
\vom

{\ubf Borel and CD cardinalities.} 
For any Borel sets $X,\,Y,$ let $X\bor Y$ mean that 
there is a Borel injection $\vt:X\to Y.$ 
Accordingly, let $X\boe Y$ mean that both $X\bor Y$ and 
$Y\bor X,$ and $X\bol Y$ will mean that $X\bor Y$ but not 
$Y\bor X.$ 
Changing ``Borel'' to ``CD'', we obtain 
${\qdr}\zd{\qde}\zd{\qdl},$ stronger relations between 
\qd\ sets. 

Obviously $X\bor Y$ iff ${\rav X}\bor{\rav Y},$ thus, the 
\dd\bor structure of Borel sets is in a sense equal to the 
\dd\bor structure of B-smooth \eqr s, and the same 
for the CD case.

\ble
\label{CB}
Let\/ $X,\,Y$ be Borel sets. 
Then\/ $X\boe Y$ iff there is a Borel bijection of\/ $X$ 
onto\/ $Y.$ 
Similarly, if\/ $X,\,Y$ are CD sets then\/ $X\qde Y$ iff 
there is a CD bijection of\/ $X$ onto\/ $Y.$ 
\ele
\bpF
Apply the Cantor -- Bernstein argument. 
To see that it yields a bijection of necessary type, recall 
that the image $\ran\vt$ of a CD, resp., Borel injection 
$\vt$ is equal to $\dom{(\vt\obr)},$ hence, is still a CD, 
resp., Borel set \cite[2.10]{kkml}.
\epF

Thus, $X\qde Y$ can be interpreted as saying that the sets 
$X,\,Y$ have the same {\it CD-cardinality\/}; the latter 
then can be defined as the \dd\qde class of $X.$ 
Similarly, $X\boe Y$ means that $X,\,Y$ have the same 
{\it Borel cardinality\/}.

The following result presents an alternative description of 
the relations ${\qde}\zd{\boe}$ restricted to $\adN$ 
(\ie, acting only on hyperfinite sets; recall that any 
$x\in\adN$ is identified with the set $[0,x)$). 

\bpr[{{\rm \cite[\S~2]{kkml}}}]
\label{cdb}
Suppose that\/ $x,\,y\in\adN.$ 
Then\/ $x\boe y$ iff\/ $\st{\frac xy}=1,$ and\/ 
$x\qde y$ iff\/ $0<\st{\frac xy}<+\infty$.\qeD
\epr

It follows that the relations $x\boe y$ and $x\qde y$ 
on $\adN$ are Borel. 
We show below that the first of them is not CD-smooth 
while the other one is CD-smooth but not B-smooth. 
\vom

{\ubf Exponential equalities.} 
Let $X,\,Y$ be nonstandard sets. 
A function $f:Y\to X$ is {\it internally extendable\/} if 
$f=g\res Y$ for an internal function $g$ with $Y\sq\dom g.$ 
(If $X$ is internal then this is the same as an internal 
function.)  

How many there are internally extendable functions $X\to2$? 
Equivalence relations allow to approach this quention 
in terms of Borel and CD reducibility.
For any nonstandard set $X,$ let $\Exp X$ be the equivalence 
relation defined on $\sui{(\fuk2H)}$ 
for some internal $H\supseteq X$ 
so that $\xi \Exp X\eta$ iff ${\xi\res X}={\eta\res X}.$ 
This definition formally depends on $H,$ but easily all 
ERs obtained this way (for a fixed $X$) are 
\dd\boe equivalent to each other, hence, 
$\Exp X$ manifests this \dd\boe type. 
If $X=H$ is itself internal then so is $\Xi=\sui{(\fuk2X)},$ 
and the definitions of $\rav{\Xi}$ and $\Exp X$ 
give obviously one and the same (modulo $\boe$).  
If $X$ is not internal then $\Exp X$ simulates the equality 
of internally extendable maps $X\to2,$ so that, for instance, 
the inequality ${\Exp X}\bor\rav Y$ ($Y$ internal) means 
that, the number of all internally extendable maps $X\to2$ 
is, in a sense, smaller-or-equal to $\#Y$.

\seci{Borel cardinalities}
\label{nB}

Our first goal is to study the \dd\bor structure of 
Borel sets in $\adN.$ 
The following theorem shows that any infinite Borel 
subset of $\adN$ is \dd\boe equivalent to a unique 
Borel cut of some kind.

\bte
\label{cardB}
For any Borel set\/ $X\sq\adN$ there is a Borel cut\/ 
$U\sq\adN$ with\/ $X\boe U,$ actually, there is a 
minimal Borel cut\/ $U$ satisfying\/ $X\boe U$.
\ete

We precede the proof of the theorem by two auxiliary
lemmas.
The first of them says that $\bor$ is sometimes
preserved under unions and intersections.

\ble[{{\rm Essentially from Zivaljevic~\cite z}}]
\label{sim}
Suppose that\/ $A_n,\,B_n$ are hyperfinite sets,
and\/ $b_n=\#B_n\le a_n=\#A_n$ for each\/ $n.$
Then
\ben
\tenu{{\rm(\roman{enumi})}}
\itla{sm1}
if\/ $A_{n+1}\sq A_n$ and\/ $B_{n+1}\sq B_n$ for each\/ $n$ 
then\/ $\bigcap_nB_n\bor\bigcap_nA_n\;;$

\itla{sm2}
if\/ $A_n\sq A_{n+1}$ and\/ $B_n\sq B_{n+1}$ for each\/ $n$ 
then\/ $\bigcup_nB_n\bor\bigcup_nA_n$.
\een
\ele
\bpF
\ref{sm1}
For any $n$ there is an internal bijection 
$f:A_0\,\text{ onto }\,[0,a_0)$
such that $f\ima A_k=[0,a_k)$ for all $k\le n.$
By \Sat, there is an internal bijection
$f:A_0\,\text{ onto }\,[0,a_0)$ with $f\ima A_n=[0,a_n)$ 
for all $n\in\dN.$
We conclude that $\bigcap_nA_n\boe U=\bigcap_n[0,a_n).$
Also, $\bigcap_nB_n\boe D=\bigcap_n[0,b_n).$
However $D\sq U.$ 

\ref{sm2}
Arguing the same way, we prove that
$\bigcup_nA_n\boe U=\bigcup_n[0,a_n)$ and 
$\bigcup_nB_n\boe D=\bigcup_n[0,b_n),$ but again $D\sq U$.
\epF

If $U\sq V\sq\adN$ are cuts then we write $U\approx V$ iff 
$\frac xy\simeq1$ for all $x,\,y\in V\bez U.$ 
(For instance, if $U=[0,a)$ and $V=[0,b)$ then $U\approx V$ 
iff $\frac ab\simeq1.$) 
This turns out to be a necessary and sufficient condition 
for $U\boe V.$

\ble
\label{U.V}
{\rm(i)} 
If\/ $U,\,V$ are Borel cuts then\/ $U\boe V$ iff\/ 
$U\approx V.$ 

{\rm(ii)} 
Any\/ \dd\approx class of Borel cuts contains a\/ 
\dd\sq minimal cut, in particular, any \poq{additive} 
Borel cut is\/ \dd\approx isolated, \ie,\/ $U\not\approx V$ 
for any cut\/ $V\ne U$.
\ele
\bpF
(i) 
Let, say, $U\sq V.$ 
Suppose that $U\boe V.$ 
Take any $x<y$ in $V\bez U.$ 
Then $x\boe y,$ hence, $\frac xy\simeq1$ by 
Proposition~\ref{cdb}. 
Suppose, conversely, that $U\approx V.$ 
Take any $x\in V\bez U.$ 
Let $c$ be the entire part of $x/2;$ then easily $c\in U.$ 
Let $A=\ans{a\in\adN:\frac ac\simeq 0}.$ 
We observe that $A\sneq U$ and the difference $D=V\bez U$ 
satisfies $D\sq X^+\cup X^-,$ where 
$X^+=\ans{x+a:a\in A}$ and $X^-=\ans{x-a:a\in A}.$ 
Define $f(z)$ for any $z\in V$ as follows. 
If $z\in U\bez A$ then $f(z)=z.$ 
If $z\in D\cap X^+$ then $z=x+a\zt a\in A,$ and we define 
$f(z)=3a$ (a number in $A$). 
If $z\in D\cap X^-,$ but $z\ne x,$ then 
$z=x-a\zt a\in A\bez\ans0,$ and we define 
$f(z)=3a+1$ (still a number in $A$). 
Finally, if $x\in A$ then let $f(x)=3x+2.$ 
Easily $f$ is a Borel injection $V\to U.$ 

(ii) 
Let $\wU$ be the set of all $x\in U$ such that 
there is $y\in U\zt y>x$ with $\frac xy\not\simeq 1.$ 
This is a cut, moreover, a projective set, hence, 
\qd, which implies that $\wU$ is actually Borel 
by Lemma~\ref{cdcut}. 
Easily $\wU\approx U.$ 
Finally, note that for any $x\in\wU$ there exists 
$x'\in\wU\zt x'>x,$ with $\frac{x'}x\not\simeq1:$ 
indeed, let $x'=\frac{x+y}2,$ where 
$y\in U\zt y>x\zt\frac yx\not\simeq 1.$ 
This suffices to infer that $V\not\approx \wU$ 
for any cut $V\sneq \wU.$ 
In other words, $\wU$ is the \dd\sq least cut 
\dd\boe equivalent to $U,$ as required.
That $\wU=U$ for any \poq{additive} cut $U$ is a simple 
exercise.
\epF 

\bpF[{{\sl Theorem~\ref{cardB}\/}}]  
Lemma~\ref{U.V} allows us to concentrate on the first 
assertion of the theorem.
Since all Borel sets are \qd, we can present a given 
Borel set $X\sq\adN$ in the   
form $X=\bigcup_{f\in F}\bigcap_nX_{f\res n},$ 
where $F$ and the sets $X_s\sq\adN$ are as in \eqref{cd} 
of Section~\ref N.
If there is $f\in F$ such that all sets
$X_{f\res n}$ are unbounded in $\adN$ then, by \Sat,
there is an {\it internal\/} unbounded set
$Y\sq X_f=\bigcap_nX_{f\res n}.$ 
Then obviously $Y\boe\adN,$ hence, $X\boe\adN$.

We assume henceforth that $X$ \poq{is bounded in} $\adN$ 
--- then it can be assumed that all sets $X_s$ are 
also bounded, hence, hyperfinite. 
Let $\nu_s=\# X_s$. 

Let $C$ be the set of all $c\in\adN$ such that there is
$f\in F$ and an internal injection
$\vpi:[0,c)\to X_f=\bigcap_nX_{f\res n}.$
Easily $C$ is a cut, and a \qd\ set.
(By \Sat, for any internal $Y$ to be internally
embeddable in $X_f$ it suffices that $\#Y\le\nu_{f\res m}$
for any $m.$)

{\it We claim that\/ $C\bor X.$\/} 
Indeed if there is $f\in F$ such that 
$C\sq[0,\nu_{f\res n})$ for all 
$n$ then immediately $C\bor X_f$ by Lemma~\ref{sim}\ref{sm1}.
Otherwise for any $f\in F$ there is $n_f\in\dN$ such that 
$\nu_{f\res n_f}\in C.$
As $X_{f\res n_f}$ is an internal set with
$\#X_{f\res n_f}=\nu_{f\res n_f},$
no internal set $Y$ with $\#Y>\nu_{f\res n_f}$ admits an
internal injection in $X_f,$
hence, the countable set $\ans{\nu_{f\res n_f}:f\in F}$
is cofinal in $C,$ so that $C=\bigcup_k[0,z_k),$ where all
$z_k$ belong to $C.$
However for any $k$ there is an internal $R_k\sq X$ with 
$\#R_k=z_k.$
Lemma~\ref{sim}\ref{sm2} implies $C\bor\bigcup_kR_k$.

In continuation of the proof of the theorem, 
we have the following cases.\vom

{\sl Case 1\/}:
$C$ is not additive. 
Then there is $c\in C$ such that $c\dN= U$ and $2c\nin C.$ 
Prove that $X\bor c\dN.$ 
By Lemma~\ref{sim}\ref{sm2}, it suffices to cover $X$ by 
a countable union $\bigcup_jY_j$ of internal sets $Y_j$ 
with $\#Y_j\le 2c$ for all $j.$
For this it suffices to prove that for any $f\in F$ there 
is $m$ such that $\nu_{f\res m}=\#X_{f\res m}\le 2c.$
To prove this, assume, on the contrary, that $f\in F$ and 
$\nu_{f\res m}\ge 2c$ for all $m;$ 
we obtain, by \Sat, an internal 
subset $Y\sq X_f$ with $\#Y=2c\nin C,$ contradiction. 
We return to this case below.\vom

In the remainder, we assume that $C$ is additive.\vom

{\sl Case 2\/}: 
$C$ is countably cofinal. 
Arguing as in Case~1, we find that for any $f\in F$ there 
is $m$ such that $\nu_{f\res m}=\#X_{f\res m}\in C.$ 
(Otherwise, using \Sat\ and the assumption of countable 
cofinality, we obtain an internal subset $Y\sq X_f$ with 
$\#Y\nin C,$ contradiction.)
Thus, $X$ can be covered by a countable union $\bigcup_jY_j$ 
of internal sets $Y_j$ with $\#Y_j\in C$ for all $j.$
It follows, by Lemma~\ref{sim}\ref{sm2}, that $X\bor C.$ 
Since $C\bor X$ has been established, we have 
$X\boe C,$ so that $U=C$ proves the theorem.\vom

{\sl Case 3\/}:
$C$ is countably coinitial, and there exists a decreasing 
sequence $\sis{h_k}{k\in\dN},$ coinitial in $\adN\bez U,$ such 
that $\frac{h_k}{h_{k-1}}$ is infinitesimal for all $k\in\dN.$
For any $k\in\dN,$ if $f\in F$ then there is
$m$ with $\nu_{f\res m}\le h_{k+1}$
(otherwise, by \Sat, $X_f$ contains an internal subset 
$Y$ with $\#Y>h_{k+1},$ contradiction),
so that $X$ is covered by a countable union of internal 
sets $Y_j$ with $\#Y_j\le h_{k+1}$ for all $j.$
It follows, by \Sat\ and because $\frac{h_k}{h_{k-1}}$ 
is infinitesimal, that, for any $k,$ 
$X$ can be covered by an internal set 
$R_k$ with $\#R_k\le h_k.$ 
Now $X\bor C$ by Lemma~\ref{sim}\ref{sm1}, hence,
$U=C$ proves the theorem.\vom

{\sl Case 4\/}:
finally, $C=c/\dN$ for some $c\nin U.$
We have $c/\dN\bor X\bor c\dN$
(similarly to Case 2).\vom

To conclude, cases 2 and 3 led us directly to the result 
required, while cases 1 and 4 can be summarized as follows: 
there is a number $c\in\adN\bez\dN$ such that 
$c/\dN\bor X\bor c\dN.$ 
We can assume that $X\sq c\dN.$ 

Let $\mu(Y)=\frac{\#Y}c$ be the counting measure on $c\dN.$ 
The set $X$ is Borel, hence, Loeb-measurable. 
If its Loeb measure is $\infty$ then there is a  
sequence $\sis{X_n}{}$ of internal subsets of $X$ with 
$\#X_n=nc\zd \kaz n.$ 
It follows that $c\dN\bor X$ by Lemma~\ref{sim}, hence, 
$X\boe U=c\dN,$ as required.

Suppose that the Loeb measure of $X$ is a (standard) 
real $r\ge0.$  
There is an increasing sequence $\sis{A_n}{n\in\dN}$ 
of internal subsets of $X$ and a decreasing sequence 
$\sis{B_n}{n\in\dN}$ of supersets of $X$ such that 
$\mu(B_n)-\mu(A_n)\to0$ as $n\to\iy$ 
(\ie, the difference is eventually less than 
any fixed standard $\vep>0$). 
If $r=0$ then $\frac{\#B_n}c\to0,$ 
therefore, $\bigcap_nB_n\boe c/\dN$ by Lemma~\ref{sim}, 
which implies $X\boe c/\dN$ since $c/\dN\bor X,$ 
therefore, $U=c/\dN$ proves the theorem.

Finally, assume that $r>0.$ 
Prove that then $X\boe [0,\ent{cr}).$ 
We have $\frac{\#A_n}c\to r$ from below and 
$\frac{\#B_n}c\to r$ from above. 
Let $U=\bigcup_{n\in\dN}[0,\#A_n)$ and 
$V=\bigcap_{n\in\dN}[0,\#B_n);$ then 
$\bigcup_nA_n\boe U$ and $\bigcap_nB_n\boe V$ by 
Lemma~\ref{sim}, while $\ent{cr}\in V\bez U,$ 
hence, in remains to prove that $U\boe V.$ 
It suffices, by Lemma~\ref{U.V},  
to show that $U\approx V.$ 
Let $x<y$ belong to $V\bez U.$ 
If $\frac yx\not\simeq1$ then $\frac yc-\frac xc$ is not 
infinitesimal, which contradicts the fact that 
$\mu(B_n)-\mu(A_n)\to0$ because $\mu(A_n)\le\frac xc$ and 
$\frac yc\le\mu(B_n)$ for all $n$.\vom

\epG{Theorem~\ref{cardB}}

\bcor
\label{B<>}
Any two Borel sets\/ $X,\,Y\sq\adN$ are\/ 
\dd\bor comparable.\qeD
\ecor

\bcor[{{\rm originally Zivaljevic~\cite z}}]
\label{Bz}
If\/ $c\in\adN\bez\dN,$ $\mu$ is a finite counting measure  
on\/ $[0,c),$ and sets\/ $X,\,Y\sq c\dN$ are Borel and 
of non-0 Loeb measure\/ $L(\mu)$ 
then\/ $X\boe Y$ iff\/ $L(\mu)(X)=L(\mu)(Y)$.
\ecor
\bpF
See the last paragraph of the proof of the theorem.
\epF

{\ubf Complete classification of Borel cardinalities.} \
Call a Borel cut $U\sq\adN$ {\it minimal\/} 
\index{cut!minimal}%
if $V\not\boe U$ for any cut $V\sneq U.$ 
It follows from Theorem~\ref{cardB} that any \dd\boe class 
of Borel subsets of $\adN$ contains a unique minimal Borel 
cut, so that minimal Borel cuts can be viewed as 
{\it Borel cardinals\/} (of Borel subsets of $\adN$). 

For instance, any additive Borel cut is minimal by 
Lemma~\ref{U.V}, hence, a Borel cardinal. 
But if $U$ is a non-additive minimal Borel cut, then 
there is a number $c\in U$ with $2c\nin U,$ so that 
$c/\dN\sneq U \sneq c\dN,$ and, accordingly, 
$c/\dN\bol U \bol c\dN,$ because $c/\dN$ and $c\dN$ are 
minimal cuts themselves. 
(Easily $c\dN$ is the least attitive cut bigger than $c/\dN.$) 

To study the structure of minimal Borel cuts between 
$c/\dN$ and $c\dN$ for a fixed nonstandard $c\in\adN,$ 
put $y_{cr}=\ent{cr}$ for any real $r\in\dR\zt z>0,$ 
where, we recall, $\ent$ is the entire part in the 
internal universe.  
Let $U_{cr}=[0,y_{cr}].$  
Easily any minimal Borel cut $U$ satisfying 
$c/\dN\bol U \bol c\dN$ is equal to $\wU_{cr}$ for 
some positive real $r,$ 
and $\wU_{cr}\ne \wU_{cr'}$ for different 
$r,\,r'$ (and one and the same $c$).
Thus, Borel cardinals of Borel subsets of 
$\adN$ are either additive Borel initial segments 
or those of the form $\wU_{cr},$ or, finally, 
(finite) natural numbers.

\seci{CD cardinalities}
\label{sCD}

It can be expected that different Borel cardinalities 
are ``glued'' by \qd\ maps. 
Lemma~\ref{carT} below reveals the exact measure of 
this phenomenon. 
The other side of the CD cardinalities vs.\ the Borel ones is 
that this notion is addressed to a much bigger class of sets, 
the \qd\ sets, which are not necessarily Loeb measurable and, 
generally, have more vague nature. 
In particular, the \dd\qdr structure of \qd\ sets is known 
only partially.

\bte
\label{tCD}
If\/ $X\sq\adN$ is an infinite \qd\ set then either 
there is a unique additive Borel cut\/ $U\qde X$ or 
there is an infinitely large\/ $c\in\adN$ such that\/ 
$c/\dN\qdl X\qdl c\dN$.
\ete

Thus, any infinite \qd\ subset of $\adN$ either is 
\dd\qde equi\-valent to a unique additive CD cut, or at least 
can be placed between two adjacent additive CD cuts. 
While the ``either'' case is realized on simple examples 
(for instance, additive Borel cuts themselves), 
the ``or'' case remains enigmatic. 

The next lemma comprises several facts involved in the proof.

\ble
\label{carT}
If\/ $U\sq\adN$ is an infinite Borel cut then\/ 
$U\qde\dN\ti U$ {\rm(the Cartesian product)} 
and\/ $U\qde U\dN$ {\rm(a cut)}. 

On the other hand, if\/ $U\sneq V$ are Borel cuts, 
and\/ $U$ is additive, 
then there is no CD map\/ $\vpi:U\,\text{ onto }\,V.$ 

It follows that, for\/ $x,\,y\in\adN,$ 
$[0,x)\qde[0,y)$ iff\/ $0<\st{\frac xy}<\iy$.
\ele
\bpF
Theorem~\ref{sq} below implies 
that there exists a CD set $W\sq U$ such that for any
$x\in U$ there is a unique $w_x\in W$ with
$|x-w_x|\in\dN.$
Let $a\mapsto \ang{z_a,n_a}$ be a recursive bijection of
$\dZ$ (the integers) onto $\dZ\ti\dN.$
Now, if $x\in U\bez\dN$ 
then put $a=x-w_x$ and $\vt(x)=\ang{w_x+z_a,n_a}.$
If $x=m\in\dN$ then let $\Phi(x)=\ang{i_m,j_m},$ where
$m\mapsto \ang{i_m,j_m}$ is a fixed bijection of
$\dN$ onto $\dN\ti\dN.$
Also, if $U$ has a maximal element $\mu$ and
$x=\mu-m\zt m\in\dN,$ then let $\vt(x)=\ang{\mu-j_m,i_m}.$
Easily $\vt$ is a CD bijection of $U$ onto $U\ti\dN$.

In the second equivalence, 
if $U$ is additive then $U=U\dN$ and there is nothing to 
prove. 
Otherwise there is $c\in U$ such that $U\dN=\bigcup_n[0,cn).$
Note that $U\dN=\bigcup_{n\in\dN}U_n,$
where $U_n=cn+U,$ hence, there is a Borel
bijection of $U\ti \dN$ onto $U\dN.$

To prove the second assertion, let, on the contrary, 
$P=\bigcup_{f\in F}\bigcap_m P_{f\res m}$
be such a map
($P_s\sq\adN\ti\adN$ are internal sets and $P_s\sq P_t$
whenever $t\subset s.$)
Then any $P_f=\bigcap_m P_{f\res m}$ is still a function, 
hence, by \Sat, there is a number $m_f$ such that 
$P_{f\res m_f}$ is a function.
Thus, there is a countable family of {\it internal\/}
functions $\Phi_i\zt i\in\dN,$ with $U\sq{\dom {\Phi_i}},$
such that $V\sq\bigcup_i\Phi_i\ima U.$
We can assume that $V=[0,c),$ where $c\in\adN\bez U.$
Put $c_0=c$ and, by induction, let $c_{n+1}$ be the 
entire part of $c_n/2.$
Then still $c_n\nin U$ for any $n$ as $U$ is an additive
cut, therefore, $V\sq\bigcup_i\Phi_i\ima[0,c_{i+2}].$
Yet every $V_i=\Phi_i\ima[0,c_{i+2})$ is an internal set
with $\#V_i\le c/2^{i+2},$ hence, by \Sat, $\bigcup_iV_i$ can
be covered by an internal set with $c/2$ elements, and cannot
cover $V.$ 
\epF

\brem
Thus, for any infinitely large $c\in\adN,$ all Borel 
cardinals (as defined in the end of Section~\ref{nB}) 
between $c/\dN$ 
and $c\dN$ are \dd\qde equivalent to each other and to 
$c\dN.$ 
It follows that for any \poq{Borel} set $X\sq\adN$ there 
is a unique {\it additive\/} Borel cut $U$ with $X\qde U,$ 
so that we can define {\it CD-cardinals\/} of 
\poq{Borel} sets to be just additive Borel cuts in $\adN.$ 
What about CD-cardinalities of \qd\ sets~? 
Unfortunately, this question remains open. 
\erem

\bpF[{{\sl Theorem~\ref{tCD}\/}}]
We leave it as an easy exercise for the reader to verify that 
the arguments in the proof of Theorem~\ref{cardB} are partially 
applicable to any \qd, not necessarily Borel, set $X\sq\adN.$ 
More exactly. 
If $X$ is unbounded in $\adN$ then $X\qde\adN.$ 
If $X$ is bounded in $\adN$ then either $X$ is 
\dd\qde equivalent to an additive Borel cut 
(cases 2 and 3) or there is an infinitely large number $c$ 
with $c/\dN\qdl X\qdl c\dN$ (cases 1 and 4). 
The Loeb measurability of Borel sets allowed us to further 
study the ``or'' case provided $X$ is a Borel set, but 
the method does not seem to apply for CD sets in general.\vom

\epG{Theorem~\ref{tCD}}

\seci{On ``singular'' CD sets}
\label{sing}

Recall that the {\it CD-cardinality\/} of a \qd\ set $X$ is 
the \dd\qde class of $X.$  
For the moment, let us consider only the case of 
{\it bounded\/} CD sets $X\sq\adN.$ 
Natural (finite) numbers and CD-cardinalities of  
additive \qd\ cuts $U\sq\adN$ can be called {\it regular\/}, 
other {\it singular\/}. 
For instance, any $c=[0,c)\in\adN$ has regular 
CD-cardinality because, by the above, if $c\nin\dN$ then 
$[0,c)\qde c\dN.$ 

\bqu
\label{?cuts}
Do there exist singular CD-cardinalities~? 
In other words (we refer to Theorem~\ref{tCD}), 
given $c\in\adN\bez\dN,$ 
does there exist a \qd\ set $X$ of type ``or'' of 
Theorem~\ref{tCD}, \ie, satisfying 
$c/\dN\qdl X\qdl c\dN$? 
If yes then are there \dd\qdr incomparable sets of this 
sort~?
\equ

If the first question answers in the negative then the 
structure of 
CD cardinalities of (\qd) subsets of $\adN$ turns out 
to be rather well organized: any infinite CD set $X\sq\adN$ 
is \dd\qde equivalent to an additive CD cut in $\adN.$ 
But we would rather conjecture the existence of ``singular'' 
\qd\ sets, \ie, those of type ``or'' of 
Theorem~\ref{tCD}. 
%
The goal of this Section is to prove that CD subsets of 
$X\sq c\dN$ satisfying $c/\dN\qdr X$ (including possible 
examples for the problem) are \dd\qde equivalent to sets of 
rather simple form, which may lead to more fruitful 
further studies.

Since any $c\in\adN\bez\dN$ belongs to an interval of the form 
$[2^d,2^{d+1})\zt d\in\adN\bez\dN,$ and then $c\dN=2^d\dN$ and 
$c/\dN=2^d/\dN,$ we can assume that already $c=2^d.$ 
In this case, the domain $c=[0,c)$ can be identified with the 
set $\Xi=\fuk2d$ of all interhal $\xi:d\to2:$ 
the map $\xi\longmapsto x(\xi)=\sum_{k=0}^{d-1}2^{d-k-1}\xi(k)$ 
is an internal bijection of $\Xi$ onto $[0,c).$ 
For any $s\in\dln$ put $M^d_s=\ans{\xi\in\Xi:s\subset\xi}.$ 
For any $g\in\dn,$ put 
$M^d_g=\bigcap_mM^d_{g\res m}=\ans{\xi\in\Xi:\xi\res\dN=g}.$ 
Call sets $M^d_g$ \dd d{\it monads\/}. 
In different terms, the monads $M^d_g$ are equivalence 
classes of the equivalence relation $\Exp\dN$ on 
$\Xi=\fuk2d,$ see Section~\ref{prel}.~\footnote
{\ In the notation of \cite{kkml}, $\xi\res\dN$ is denoted 
by $\st\xi,$ the standard part, hence, we have 
$M^d_g=\st\obr(\ans g)$ and $M^d_G=\st\obr(G).$}

For instance, $M^d_{\bf 0},$ where ${\bf 0}\in\dn$ is the 
constant $0,$ is a \dd dmonad. 
Easily $\ans{x(\xi):\xi\in M^d_{\bf 0}}=c/\dN,$ hence, 
${c/\dN}\qde{M^d_{\bf 0}}\qde{M^d_g}$ for each $g\in\dn.$

Any union $M^d_G=\bigcup_{g\in G}M^d_g$ of \dd dmonads 
($G\sq\dn$) is clearly a CD set.

\bte
\label{rer-c}
Suppose that\/ $c=2^d\in\adN\bez\dN.$ 
If\/ $X\sq c\/\dN$ is a \qd\ set then either\/ 
$X\qdl c/\dN$ or\/ $X\qde M^d_G$ for some\/ $G\sq\dn$.
\ete
\bpF
As $c=[0,c)\qde c\dN$ by Lemma~\ref{carT}, we can assume that 
$X\sq c,$ moreover, $X=\bigcup_{f\in F}\bigcap_mX_{f\res m},$ 
where $F\sq\dn$ while $X_s\sq c$ are internal sets. 
We claim that the following can be \noo\ assumed:
\ben
\tenu{(\arabic{enumi})}
\itla{ass1}\msur 
$X_t\sq X_s$ whenever $s\subset t$ 
(otherwise put $X'_s=\bigcap_{k<\lh u}X_{s\res k}$); 

\itla{ass2}\msur  
$X_{s\we0}\cap X_{s\we1}=\pu$ 
for any $s\in\dln;$

\itla{ass3} 
for any $s\in\dln,$ either $\#X_s=c2^{-k}$ for some $k=k_s$ 
or $\#X_s\in c/\dN$;

\itla{ass4} 
for any $s\we i\in\dln,$ $\#X_{s\we i}\le\frac12\#X_s$ 
($s\we i$ is the extension of $s$ by $i$). 
\een

{\sl Justification of \ref{ass2}\/}.
Sets $X_s$ admit partitions 
$X_s=\bigcup\cX_s,$ where $\cX_s$ is a finite collection of 
pairwise disjoint internal subsets of $X_s$ such that 
\ben
\tenu{(\alph{enumi})}
\itla{ben1}
if $s\subset t$ then for any $A\in \cX_t$ there is (unique) 
$B\in \cX_s$ with $A\sq B$;

\itla{ben2}
if $s,\,t\in\dln$ have the same length then any $A\in\cX_s$ 
and $B\in\cX_t$ are either equal or disjoint.
\een
Now, let $\Phi$ be the set of all functions $\vpi\zt\dom\vpi=\dN,$ 
such that there is $f\in F$ satisfying 
$\vpi(m)\in\cX_{f\res m}$ and 
$\vpi(m+1)\sq\vpi(m)$ for any $m.$ 
Obviously $X=\bigcup_{\vpi\in\Phi}\bigcap_m\vpi(m),$ which  
justifies the claim by a modification of the argument used to 
derive \eqref{cd} from \eqref{cd1} in Section \ref N. 

{\sl Justification of \ref{ass3}\/}.
Partitions $X_s=\bigcup\cX'_s$ can be defined, such that  
$\cX'_s$ is an at most countable collection of subsets of 
$X_s,$ of which at most one, say $P_s,$ is a $\Pi^0_1$ set 
with $P_s \qdr c/\dN$ while 
all others are (pairwise disjoint) internal sets of 
hyperfinite cardinalities of the form $c2^{-k}\zt k\in\dN,$ 
and still both \ref{ben1} and \ref{ben2} hold 
(for the collections $\cX_s'$). 
We can drop all sets $P_s$ because this amounts to a total set  
of $\qdr c/\dN$ elements by Lemma~\ref{carT}, which is not 
essential in the context of the theorem. 
Then proceed as above. 

{\sl Justification of \ref{ass4}\/}: 
a similar argument.

Coming back to the proof of the theorem, let 
$S=\ans{f\res m:f\in F\land m\in\dN}$ (a subset of $\dln$). 
In the assumptions \ref{ass1} -- \ref{ass4}, one can define 
$\sg_s\in\dln$ for any $s\in S$ so that 
(A) if $s\we i\in S$ ($i=0\text{ or }1$) then 
$\sg_s\we i\sq \sg_{s\we i},$ 
and 
(B)~if $\#X_s=c\/2^{-m}$ then $\lh {\sg_s}=m.$ 
For any $f\in F,$ let $g(f)=\bigcup_m\sg_{f\res m}\in\dn.$ 
%
Let $G=\ans{g(f):f\in F}.$ 
Note that (in our assumptions) the sets 
$X_f=\bigcap_mX_{f\res m}$ and 
$M^d_{g(f)}=\bigcap_nM^d_{g(f)\res n}=\bigcap_mM^d_{\sg_{f\res m}}$ 
are $\qde$ by Lemma~\ref{sim}\ref{sm1}, moreover, by a 
suitable modification of the proof of Lemma~\ref{sim}\ref{sm1}, 
we find an internal map $\vt$ such that 
$\vt\ima X_f=M^d_{g(f)}$ for any $f\in F,$ 
hence, $\vt\res X$ is a bijection of $X$ onto $M^d_G.$  
\epF

Thus, if Problem~\ref{?cuts} answers in the positive 
then, by the theorem, there exist corresponding examples 
of the form $M^d_G\zt G\sq\dn.$ 
The following rather elementary consideration focuses on 
\dd\qdr properties of sets of this form.
 
Let a number $c=2^d\in\adN\bez\dN$ be fixed. 

First of all note that $M^d_g\qde c/\dN$ for any $g\in\dn,$ 
see above, therefore, we have $c/\dN\qdr M^d_G$ whenever 
$\pu\ne G\sq\dn$.  

For $c/\dN\qdl M^d_G$ (strictly), 
it is necessary and sufficient that $\leB G>0,$ where $\leB$ 
is the upper Lebesque measure on $\dn.$ 
(Indeed, if $\leB G=0$ then, by \Sat, for any $m$ the set 
$M^d_G$ can be covered by an internal set $X$ with $\#X\le c/m,$ 
hence, $c/\dN\qdl M^d_G$ by Lemma~\ref{sim}\ref{sm1}. 
Conversely, if $M^d_G$ is covered by an internal set $X$ with 
$\#X\le c/m,$ then, for any $g\in G,$ there is a number 
$m_g\in\dN$ such that 
$M^d_{g'}\sq X$ whenever $g'\in D_{g\res m_g},$ 
where $D_s=\ans{g'\in\dn:s\subset g'}$ for any $s\in\dln.$ 
But the union $D=\bigcup_{g\in G} D_{g\res m_g}$ easily has 
measure $\le m^{-1}$ in $\dn.$)

If $\leb G>0$ (the inner Lebesgue measure) then 
$M^d_G\qde[0,c)\qde c\dN.$ 
Indeed, we can assume, by Cantor -- Bernstein, that $G$ is a 
closed subset of $\dn$ of positive measure, say, of measure 
$2^{-m}\zt m\in\dN.$ 
Then $M^d_G$ is equal to a decreasing union $\bigcap_nX_n,$ 
where each $X_n$ is an internal set with $\#X_n\ge c2^{-m}.$ 
It follows, by Lemma~\ref{sim}, that $[0,c2^{-m})\qdr M^d_G,$ 
and so on.

But $\leb G>0$ is not a necessary condition 
for $M^d_G\qde c\dN.$
Indeed, let $G$ be a transversal for the equivalence relation 
$f\Eo g$ iff $f(n)=g(n)$ for all but finite $n\in\dN$   
($f,\,g\in\dn$), 
an example of a set with $\leb G=0$ and $\leB G=1.$ 
There is a sequence of internal functions $\vt_n$ such 
that $[0,c)=\bigcup_n\vt_n\ima M^d_G,$ so that, by an 
argument similar to Lemma~\ref{carT}, we have 
$M^d_G\qde [0,c)\qde c\dN.$ 

Thus, to obtain an anticipated example for 
Problem~\ref{?cuts} in the form $M^d_G,$ we have to employ 
nonmeasurable sets $G\sq\dn$ with $\leb G=0<\leB G$ 
but less ``dense'' than transversals of $\Eo.$ 
It remains to be seen whether such an approach may lead 
to a solution of the problem.

\bqu
Which ``standard'' property of $G,\,G'\sq\dn$ is 
necessary and sufficient for $M^d_G\qde M^d_{G'}$?
\equ

\seci{Countable \er s have transversals}
\label{ctble}

An equivalence relation $\qE$ is {\it ``countable''\/} 
if any its equivalence class, \ie, a set of the form 
$\ek x=\ans{y:x\qE y}\zt x\in\dom\qE,$ is at most countable. 
In \pol\ descriptive set theory, ``countable'' Borel \er s 
form a rather rich class whose full structure 
in terms of Borel reducibility is a topic of deep 
investigations (see Kechris~\cite{ndir}). 
In nonstandard setting, the picture is different.  

\bte
\label{sq}
Any ``countable'' \qd\ \eqr\/ $\qE$ on\/
$\adN$ admits a \qd\ transversal, hence, is CD-smooth.
\ete

Jin~\cite j proved the result for the \er\ 
$x\qen y$ iff $|x-y|\in\dN.$ 
Our proof of the general result employs a somewhat 
different idea, although some affinities with Jin's 
arguments can be traced. 
Note also that $\qen,$ a typical countable \eqr, is 
\poq{not} B-smooth (see Lemma~\ref{NB} below), this is 
the most transparent case when the Borel reducibility 
is really stronger.

\bpF 
The CD-smoothness easily follows from the existence of 
a transversal: just let $\vt(x)$ to be the only element 
of a transversal equivalelent to $x$.

To define a transversal, suppose, as usual,  that 
$\qE=\bigcup_{f\in F}\bigcap_{m\in\dN}P_{f\res m},$ where all
sets $P_s\zd s\in \dln,$ are internal subsets of 
$\adN\ti\adN$ with $P_t\sq P_s$ whenever $s\subset t,$ 
and $F\sq\dn.$ 
An ordinary \Sat\ argument   
shows that, because all \dd\qE classes are countable and 
a countable set cannot contain an infinite internal subset, 
for any $f\in F$ there is a number $m_f\in\dN$ such that 
all cross-sections 
$P_{f\res m_f}(x)=\ans{y:\otn x{P_{f\res m_f}}y}$ 
are finite. 
Let $S=\ans{f\res m_f:f\in F};$ this is a subset of $\dln.$ 
Then, for any $s\in S\zt k\in\dN,$ and $x\in\adN,$ we can 
define $f_{sk}(x)$ to be the \dd kth element 
(the counting begins with $0$) of $P_s(x),$ 
in the natural order of $\adN,$ whenever $\#P_s(x)\ge k,$ 
so that $f_{sk}$ is an internal partial function 
$\adN\to\adN$.

Let $s\in S$ and $k\in\dN\zt k\ge1.$ 
For any $x\in\adN$ define an internal decreasing sequence 
$\sis{x_{(a)}}{a\le a(x)}$ of length $a(x)+1\in\adN$ as 
follows. 
Put $x_{(0)}=x.$ 
Suppose that $x_{(a)}$ is defined. 
If $z=f_{sk}(x_{(a)})$ is defined and $z<x_{(a)}$ then 
put $x_{(a+1)}=z,$ otherwise put $a(x)=a$ and end the 
construction. 
(Note that eventually the construction stops simply because 
$x_{(a+1)}<x_{(a)}.$)
Put $\nu_{sk}(x)=0$ if $a(x)$ is even and $\nu_{sk}(x)=1$ 
otherwise.

Put 
$\pro x=\ans{\ang{s,k}\in S\ti\dN:\nu_{sk}(x)=0},$ 
the ``profile'' of any $x\in\adN$.

\ble
\label{ple}
If\/ $x\ne y\in\adN$ and\/ $x\qE y$ then\/ $\pro x\ne \pro y$.
\ele

Thus, while it is, generally speaking, possible that different
nonstandard numbers have equal ``profiles'', this cannot happen
if they are \dd\qE equivalent.

\bpF
We can assume that $y<x.$ 
There is $f\in F$ such that 
$\ang{x,y}\in P_f=\bigcap_mP_{f\res m}.$ 
let $s=f\res m_f,$ an element of $S.$ 
Then $y$ belongs to $P_s(x),$ a finite set, say, $y$ is 
\dd kth element of $P_s(x),$ in the natural order of $\adN.$ 
In other words, $y=x_{(1)},$ in the sense of the construction 
above, therefore, $y{(1)}=x_{(2)},$ \etc.; 
we conclude that $\nu_{sk}(x)\ne\nu_{sk}(y)$.
\epG{Lemma}

Coming back to the theorem, choose  
an element $r_A\in A$ in any set $\pu\ne A\sq\ps(S\ti\dN).$
For any $x\in\adN,$ the set $A(x)=\ans{\pro y:y\in\ek x}$
is a non-empty countable subset of $\ps(\dln\ti\dN).$
Then $X=\ans{x\in\adN:\pro x=r_{A(x)}}$ 
is a transversal for $\qE$ by Lemma~\ref{ple}.
To prove that $X$ is \qd\ consider the family $S$ which 
consists of all sets 
\dm
D_{sk}=\dom f_{sk}\,,\quad  
X_{sk}=\ans{x\in\adN:\nu_{sk}(x)=0}\,,
\dm 
and 
$X_{sks'k'}=\ans{x\in D_{sk}:\nu_{s'k'}(f_{sk}(x))=0},$ 
along with their complements. 
Let $\cA$ be the set of all at most countable sets 
$A\sq\ps(S\ti\dN).$ 
Obviously $X=\bigcup_{A\in\cA}X(A),$ where 
\dm
X(A)=\ans{x\in X:A(x)=A}
=\ans{x\in \adN:A(x)=A\land \pro x=r_A}\,,
\dm

\ble
\label{b:allcd}
Any set\/ $X(A)\zt A\in\cA,$ is \qd\ in\/ $S,$ in the sense  
that it can be obtained by\/ \eqref{cd1} of Section~\ref N
applied to sets in\/ $S.$
\ele
\bpF
Direct straightforward reduction to sets in $S$ shows that 
$X(A)$ is even Borel in $S$ in a similar sense. 
The most essential part of the reduction is to express the 
inclusion $A(x)\sq A$ by the formula
\dm
\kaz s\in S\;\kaz k\in\dN\;
\skl
x\in D_{sk}\imp \sus r\in A\:(r=\pro f_{sk}(x))
\skp\,,
\dm
to avoid a universal quantifier over the equivalence 
class $\ek x.$
\epF

On the other hand, the class of all sets \qd\ in a fixed 
{\it countable\/} 
collection $S$ of internal sets is closed under any unions 
(as well as under complements and intersections): 
just take the set theoretic union of the ``bases'' $B$   
in the assumption that the assignment of sets in $S$ to 
indices is fixed once and for all.
(Note that the class of all CD sets is closed only under
countable unions and intersections~!)
\epG{Theorem~\ref{sq}}

\bcor
\label{sq1}
The equivalence relation\/ $x\qde y$ on\/ $\adN$ 
admits a \qd\ transversal.
\ecor
\bpF
Recall that $x\qde y$ iff $0<\st{\frac xy}<+\iy,$ 
Proposition~\ref{cdb}.  
It follows that 
the set $\ans{2^x:x\in X},$ where $X$ is any CD transversal 
for the countable relation $x\qE_\dN y$ iff $|x-y|\in\dN,$ is
as required.
\epF

On the contrary, the relation $x\boe y$ 
iff\/ $\st\frac xy=1$ does not have a CD transversal.
Indeed, suppose that $X$ is a CD transversal for $\boe$
restricted to the set $D=[c,2c],$ where $c$ is a fixed 
infinitely large hyperinteger.
Note that, for $x,\,y\in D,$ $x\boe y$ is equivalent
to $\st\frac xc=\st\frac yc$ 
so that $X$ yields a CD
transversal for the equivalence relation of 
``having the same standard part $\st r$'' 
on the set of hyperrationals 
$A=\ans{r=\frac xz:x\in D},$ 
known to be impossible \cite[2.6]{kkml}. 
In fact ``the same standard part'' \er\ is not CD-smooth 
and even not \dd\qdr reducible to any $\Sg^0_1$ \er; 
this can be derived from our result in Part 2 of 
Section~\ref{ccoi}.

\seci{Silver -- Burgess dichotomy}
\label{silv}

It is known from Henson~\cite{he}
(see also Proposition 2.5 in \cite{kkml})
that any \qd\ {\it set\/} $X\sq\adN$ is countable
or else contains an infinite internal subset.
The following is a slight generalization.

\ble
\label{H}
Let\/ $X\sq\adN$ be a \qd\ set and\/ $U\sq\adN$ an additive 
cut of countable cofinality.
\ben
\tenu{{\rm(\roman{enumi})}}
\itla{H1}
Either\/ $X\qdr U$ or\/ $X$ contains an internal subset\/ $Y$ 
with\/ $\#Y\nin U$. 

\itla{H2}
Either\/ $X$ is bounded\/
{\rm(\ie, $X\sq h$ for some $h\in\adN$)}
or\/ $X$ contains an unbounded internal subset.
\een
\ele
\bpF~\footnote
{\ The result follows from Theorem~\ref{tCD}, but we prefer 
to present a short direct proof.} 
\ref{H1} 
Suppose that $X=\bigcup_{f\in F}\bigcap_n X_{f\res n},$ where 
$F\sq\dn$ and $X_s$ are as in \eqref{cd} of Section~\ref N.  
Let $S$ consist of all $s\in\dln$ with $\# X_s\in U.$   
If there is $f\in F$ such that $f\res n\nin S$ for all $n$ 
then by \Sat\ $\bigcap_n X_{f\res n}$ contains an 
internal subset $Y$ with $\#Y\nin U.$   
Otherwise we have the ``either'' case.

\ref{H2}
A similar argument, with $S$ being the set of all $s\in\dln$
such that $X_s$ is unbounded in $\adN$.
\epF

Quotient structures $\adN/{\qE},$ where $\qE$ is a
CD equivalence relation, normally consist of non-internal
elements, hence, do not contain internal subsets, but we
can consider internal pairwise \dd\qE inequivalent sets 
(\ie, sets of pairwise \dd\qE inequivalent elements) instead. 
This leads us to the following theorem, saying that, 
given a \qd\ ER $\qE,$ either the number of equivalence 
classes is somehow restricted or there is a rather big 
pairwise inequivalent set.
Recall that the relation $\Exp U$ of equality of internally 
extendable maps $U\to2$ was defined in Section~\ref{prel}.  

\bte
\label{Ts}
Let\/ $\qE$ be a CD equivalence relation on\/ $\adN,$ and\/ 
$U$ a countably cofinal additive cut. 
Then either\/ $\qE\qdr {\Exp U}$ or there is an internal 
pairwise\/ \dd\qE inequivalent set\/ $Y\sq\adN$ 
with\/ $\# Y\nin U$. 

In particular, either\/ $\qE\qdr {\Exp\dN}$ 
{\rm(then $\qE$ has \dd{{\le}\cont}many equivalence 
classes)} 
or there is an infinite internal pairwise\/ 
\dd\qE inequivalent set\/ $Y\sq\adN$.
\ete

It is not clear whether the general case is really 
a dichotomy and the quantitative characteristics are optimal. 
Generally, if $U=\bigcup_n[0,a_n)$ is countably cofinal 
($\sis{a_n}{n\in\dN}$ increases) 
then the \eqr\ $\Exp U$ has exactly 
$\prod_{n\in\dN} \card{({\pv{[a_n,a_{n+1})}})}$ 
many equivalence clas\-ses,
which is equal to $\kpa^\om$ 
provided all infinite internal sets have the same 
cardinality $\kpa$ and the differences $a_{n+1}-a_n$ are 
infinite.

The relation $\Exp\dN$ 
(the particular case in the theorem) has exactly 
\dd\cont many equivalence classes and does not
admit an infinite internal pairwise inequivalent set 
(see \cite[2.6]{kkml} on the last claim), hence, 
the continuum cannot be improved to any smaller cardinal 
in the particular case. 
 
\bpF[{{\rm Theorem~\ref{Ts}\/}}]
Suppose that 
${\qE}=\bigcup_{f\in F}\bigcap_{m\in\dN}P_{f\res m},$ 
where $P_s$ are internal subsets of 
$\adN\ti\adN$ with $P_t\sq P_s$ whenever $s\subset t,$ 
as in \eqref{cd} of Section~\ref N, while $F\sq\dn.$ 
We can \noo\ assume that the sets $P_s$ are symmetric, 
\ie, $P_s=P_s{}\obr:$ indeed, if this is not the case, then, 
as $\qE$ itself is symmetric, 
\dm
\TS
\qE={\qE}\cup{\qE\obr}=
\bigcup_{f\in F}\bigcap_{m\in\dN}
(P_{f\res m}\cup P_{f\res m}{}\obr)\,,
\dm
where the sets $P_{f\res m}\cup P_{f\res m}{}\obr$ 
are symmetric.

Since $\qE$ is an \eqr, we have
\dm
\TS
\sus f\in F\;\sus z\;\kaz m\;
(\otn x{P_{f\res m}}z \land \otn y{P_{f\res m}}z )
\,\imp \,
{\otn x\qE y}\,.
\dm
%
By \Sat, this can be rewritten as 
\envur{
\label{F2}
\TS
\kaz T\in A(F)\;\sus s\in T\;\sus z\; 
(\otn x{P_{s}}z \land \otn y{P_{s}}z )
\,\imp \,
{\otn x\qE y}\,,
}
where $A(F)$ is the collection  of all sets 
$T\sq\ans{f\res m:f\in F\land m\in\dN}$ such that 
$T\cap\ans{f\res m:m\in\dN}\ne\pu$ for each $f\in F$. 

Now let $\sis{a_n}{n\in\dN}$ be an increasing sequence 
cofinal in $U.$ 
Suppose that there is no internal pairwise 
\dd\qE inequivalent set $Y$ with $\#Y\nin U,$ more formally,
\dm
\kaz Y\in\pv\adN\:
\skl
\kaz k\:(\#Y>a_k)\imp \sus x\ne y \in Y \: 
\sus f\in F \: 
\kaz m \: (\otn x{P_{f\res m}}y)
\skp
\dm
where $\pv\adN=\ans{Y\sq\adN:Y\,\text{is internal}}.$ 
The expression to the right of $\imp$ can be 
consecutively transformed 
(using \Sat\ and the fact that $P_t\sq P_s$ 
provided $s\subset t$) 
to 
$\sus f\in F \: \kaz m\:\sus x\ne y \in Y \: 
(\otn x{P_{f\res m}}y),$ 
and then to 
\dm 
\kaz T\in A(F)\:\sus s\in T \: 
\sus x\ne y \in Y \: 
(\otn x{P_s}y)\,, 
\dm
which leads us to the following, for every $T\in A(F)$: 
\dm
\bay{l}
\kaz Y\in\pv\adN\:
\skl
\kaz k\:(\#Y>a_k)\imp \sus s\in T \: 
\sus x\ne y \in Y \: 
(\otn x{P_s}y) 
\skp\,.
\eay
\dm
Applying \Sat\ once again, we obtain, for any set $T\in A(F),$ 
a number $k(T)\in\dN$ and a finite set $T'\sq T$ 
such that 
\dm
\kaz Y\in\pv\adN\:
\skl
\#Y>a_{k(T)}\imp \sus s\in T' \: 
\sus x\ne y \in Y \: 
(\otn x{P_s}y) 
\skp\,.
\dm
Since the sets $P_s$ are assumed to be symmetric, we 
conclude that for any $T\in A(F)$ there exists an internal 
set $Z_T\sq\adN$ 
satisfying $\#Z_T\le a_{k(T)}$ and 
\envur{
\label{pop3}
\kaz x\in\adN\;\sus z\in Z_T\;\sus s\in T' \; 
(\otn x{P_s}z) \,.
}
Yet \eqref{pop3}, as a property of $Z_T,$ 
depends only on $T',$ a finite subset of $\dln,$ not 
on $T$ itself, hence, we can choose sets $Z_T$ so that 
there are only countably many different among them.  
As $U$ is an additive cut, the set 
$Z=\bigcup_{T\in A(F)}Z_T\sq\adN$ admits, by \Sat, 
an internally extendable 
(see Section~\ref{prel}) injection $\vpi:Z\to U,$ 
moreover, the cartesian product $Z\ti\dln$ also admits an 
internally extendable injection in $U,$ hence, it suffices 
to prove that $\qE\qdr{\Exp{Z\ti\dln}}$.

Let $H$ be any internal set with $Z\ti\dln\sq H.$ 
Put, for any $x\in\adN,$ 
$\vt_x=\ans{\ang{z,s}\in H:{\otn x{P_s} z}},$ 
where $P_s$ is uniformly defined via an arbitrary internal 
extension of the external map $s\mapsto P_s$ defined on 
$\dln.$ 
We have to show that, for $x,\,y\in\adN,$ 
$\vt_x=\vt_y$ implies $x\qE y.$  
Assuming that $\vt_x=\vt_y,$ fix $T\in A(F).$ 
Choose, by \eqref{pop3}, $z\in Z$ and 
$s\in T$ with $\otn x{P_s}z$ --- 
then $\ang{z,s}\in\vt_x=\vt_y,$ hence, we also have 
$\otn y{P_s}z.$ 
It remains to refer to \eqref{F2}.
\epF

Equivalence relations of class $\Sg^0_1$ admit the following 
special result:

\ble
\label{Lcof}
Assume that\/ $\qE$ is a\/ $\Sg^0_1$ \eqr\ on a subset 
of\/ $\adN,$ and\/ $X\sq\dom\qE.$ 
Then$:$
\ben
\tenu{{\rm(\roman{enumi})}}
\itla{lcof1}
if\/ $X$ is\/ $\Pi^0_1$ then either the quotient\/ 
$X/\qE$ is finite or there is an infinite internal 
pairwise\/ \dd{\qE}inequivalent set\/ $C\sq X;$

\itla{lcof2}
if\/ $X$ is \qd\ then either\/ 
$X/\qE$ is at most countable or there is an infinite 
pairwise\/ \dd{\qE}inequivalent internal set\/ $C\sq X.$
\een
\ele
\bpF
\ref{lcof1}
Let $X=\bigcap_nX_n$ and $\qE=\bigcup_nE_n,$ all $X_n$ 
and $E_n$ being internal and 
$X_{n+1}\sq X_n,$ $E_{n}\sq E_{n+1}$ for all $n.$ 
If $X/\qE$ is infinite then, for any $n,$ there is an 
internal set $C\sq X_n$ with $\#C\ge n,$ such that 
$\ang{x,y}\nin E_n$ for any two elements $x\ne y$ of $C.$ 
It remains to apply \Sat.

\ref{lcof2}
Let $X=\bigcup_{f\in F}\bigcap_m X_{f\res m},$ where 
$F$ and $X_s$ are as in \eqref{cd} of Section~\ref N. 
If for any $f\in F$ there is a number $m_f$ such that 
$X_{f\res m_f}/\qE$ is finite then $X/\qE$ is at most 
countable. 
Otherwise there is $f\in F$ such that $X_{f\res m}/\qE$ 
is infinite for all $m,$ and, arguing as in \ref{lcof1}, 
we obtain an infinite pairwise \dd{\qE}inequivalent 
internal subset of $X_f=\bigcap_m X_{f\res m}$.
\epF

\seci{Monadic equivalence relations}
\label{mer}

Any additive cut $U\sq\adN$ defines a {\ubf monadic\/} 
\eqr\ $x\qeu y$ iff $|x-y|\in U$ on $\adN.$
(If $U$ is not additive then $\qeu$ may not be a \er.)
Classes of \dd\qeu equivalence, that is, sets of the
form $[x]_U=\ans{y:|x-y|\in U}\zt x\in\adN,$
are called \dd U{\ubf monads\/},
all of them are convex subsets of $\adN.$ 

It follows from Lemma~\ref{cdcut} that there are 
two types of \qd\ monadic \er s $\qeu:$ 
{\it countably cofinal\/} and 
{\it countably coinitial\/}, 
\pagebreak[0]%
according to the type of the cut $U.$ 
(The only exceptions are $\qei\pu,$ the   
equality on $\adN,$ and $\qei\adN,$ the relation   
which makes all elements of $\adN$ equivalent.)
It turns out that the relations between monadic \er s 
in terms of $\qdr$ are determined by the relative rate 
of growth or decrease of corresponding cofinal or  
coinitial sequences. 

To distinct cuts of lowest possible rate, say that  
additive countably cofinal cuts of the form 
$c\dN\zt c\in\adN$ 
and countably coinitial cuts of the form 
$c/\dN\zt c\in\adN\bez\dN)$ are {\it slow\/},  
while other additive countably cofinal or coinitial 
cuts are {\it fast\/}.  
For instance, $\dN$ is a slow cut. 
The following is easy:

\ble
\label{?fast}
A countably cofinal additive cut\/ $U$ is slow iff\/ 
$U=\cud{r+n}$ for some\/ 
$r\in\adN,$ and is fast iff there 
is an increasing sequence\/ $\sis{a_n}{n\in\dN}$ in\/ 
$\adN$ such that\/ $U=\cud{a_n}$ and\/ 
$a_{n+1}-a_n$ infinite for any\/ $n.$ 

A countably coinitial additive cut\/ $U$ is slow iff\/ 
$U=\cud{r-n}$ for some\/ 
$r\in\adN\bez\dN,$ and is fast iff there 
is a decreasing sequence\/ $\sis{a_n}{}$ in\/ $\adN$ such 
that\/ $U=\cvd{a_n}$ and\/ 
$a_{n}-a_{n+1}$ infinite for any\/ $n.$\qeD
\ele

(Recall that $\cu{p_n}=\bigcup_{n\in\dN}[0,p_n)$ and 
$\cv{p_n}=\bigcap_{n\in\dN}[0,p_n).$)

Fast cuts admit further analysis.
If $\sis{a_n}{}$ and $\sis{b_k}{}$ are {\ubf increasing} 
sequences of hyperintegers, then define 
$\sis{a_n}{}\lek\sis{b_k}{}$ 
(meaning: $\sis{b_k}{}$ increases faster) iff 
\envur{
\label{ke+}
\kaz k\;\sus n\;\kaz n'>n\;\sus k'>k\;
\left(\frakk{a_{n'}}{a_n}\le\frakk{b_{k'}}{b_k}\right).
}
Note that the negation of \eqref{ke+} has the form
\envur{
\label{nke+}
\tag{$\mtho\neg\:\arabic{equation}$}
\sus k\;\kaz n\;\sus n'>n\;\kaz k'>k\;
\left(\frakk{a_{n'}}{a_n}>\frakk{b_{k'}}{b_k}\right).
}
Accordingly, if $\sis{a_n}{}$ and $\sis{b_k}{}$ are 
{\ubf decreasing} 
sequences of hyperintegers, then we define 
$\sis{b_k}{}\lek\sis{a_n}{}$ 
(meaning: $\sis{b_k}{}$ decreases faster) if and only if 
\envur{
\label{ke-}
\kaz k\;\sus n\;\kaz n'>n\;\sus k'>k\;
\left(\frakk{a_{n}}{a_{n'}}\le\frakk{b_{k}}{b_{k'}}\right).
}

Finally, if $U,\,V$ are countably cofinal additive cuts, 
then $U\lek V$ means that there are increasing sequences 
$\sis{a_n}{}\lek\sis{b_k}{}$ with $U=\cu{2^{a_n}}\zt 
V=\cu{2^{b_k}}.$ 
Similarly, if $U,\,V$ are countably coinitial additive cuts, 
then $U\lek V$ means that there are decreasing sequences 
$\sis{a_n}{}\lek\sis{b_k}{}$ with $U=\cv{2^{a_n}}\zt 
V=\cv{2^{b_k}}.$ 

\brem
\label{s.f}
If $U,\,V$ are both countably cofinal or both countably 
coinitial additive cuts, and $U$ is slow, 
then $U\lek V$ by Lemma~\ref{?fast}.
\erem

\bte
\label T
Suppose that\/ $U,\,V$ are additive \qd\ cuts in\/ 
$\adN$ other than\/ $\pu$ and\/ $\adN.$ 
Then\/ ${\rav\adN}\bor\qeu.$
In addition,
\ben
\tenu{{\rm(\roman{enumi})}}
\itla{T2}
if both\/ $U,\,V$ are countably cofinal or both are 
countably coinitial then\/ 
$\qeu$ and\/ $\qev$ are\/ \dd\bor comparable, and\/ 
$\qeu\bor\qev$ iff\/ $\qeu\qdr\qev$ iff\/ $U\lek V,$ 
in particular, if\/ $U$ is slow then\/ $\qeu\qdr\qev$ 
{\rm(see Remark~\ref{s.f})}$;$

\itla{T4}
if\/ $U$ is countably cofinal and\/ $V$ countably 
coinitial then\/ $\qev\not\qdr\qeu$ and\/ 
$\qeu\not\bor\qev,$ while\/ 
$\qeu\qdr\qev$ holds iff\/ $U$ is slow$;$ 

\itla{T1}\msur
$\qeu$ is not B-smooth, and\/ $\qeu$ is CD-smooth 
if and only if\/ $U$ 
is countably cofinal and slow$;$

\itla{T5}
for any countable sequence of countably cofinal  
fast cuts\/ $U_n$ there are countably cofinal   
fast cuts\/ $U,\,V$ with\/ 
$\qeu\bol\qei{U_n}\bol\qev\zd\kaz n,$ 
and the same for countably coinitial cuts.
\een
\ete

This theorem, which explains the \dd\qdr structure of 
mo\-na\-dic \eqr s, 
will be the focal point in the remainder.
According to the theorem, 
\qd\ monadic \er s form two distinct linearly 
\dd\qdr(pre) ordered domains, one of which contains 
countably cofinal and the other one countably coinitial 
\er s, each has slow \er s as the \dd\qdr least 
element, and there is no \dd\qdr connection between them 
except that any slow countably cofinal \er\ 
(it is necessarily CD-smooth) 
is \dd\qdr reducible to any countably coinitial \er. 
In addition, each of the domains is neither countably 
\dd\qdr cofinal nor countably \dd\qdr coinitial in its 
fast part. 
(It can be shown that each of the domains is also dense 
and countably saturated, \ie, contains no gaps of countable 
character.)

The proof begins with a couple of auxiliary results.

\seci{Two preliminary facts}
\label{pre}

The first result will be a connection between monadic 
\er s and certain natural \eqr s on dyadic sequences. 
Let $\fse$ be the (internal) set of all internal  
sequences $\vpi\in \ddn$ such  
that the set $\ans{a:\vpi(a)=1}$ is hyperfinite. 

Consider an additive cut $\pu\ne U\ne\adN.$ 
Then $\log U=\ans{a\in\adN:2^a\in U}$ is still 
a cut ( not necessarily additive).
Define the \eqr\ $\reil U$ on $\fse$ as 
follows: $\vpi \reil U\psi$ iff 
$\vpi\res{(\adN\bez\log U)}=\psi\res{(\adN\bez\log U)}.$ 
The relation $\reil U$ can be viewed as the restriction of 
$\Exp{\adN\bez\log U}$ (Section~\ref{prel})~to~$\fse.$ 

\bpr
\label{pcof}
In this case, $\qeu\boe\reil U$.
\epr
\bpF
For any $x\in\adN$ there is a unique $\sg=\bine x\in\fse$ with  
$x=\sum_{z\in\adN}2^z\sg(z)$ in $\adN.$ 
(The essential domain of summability here is a hyperfinite 
set because $\sg\in\fse.$)
The map $x\mto\bine x$ is not yet a reduction of 
$\qeu$ to $\reil U$ because of a little discrepancy. 
Let $\fse_{\log U}$ be the set of all $\sg\in\fse$ which are 
not eventually $1$ in $\log U,$ \ie, the set 
$\ans{a\in\log U:\sg(a)=0}$ is cofinal in $\log U.$ 
Let $\Om_{\log U}$ be the set of all $x\in \adN$ such that
$\bine x\in\fse_{\log U}.$ 

We assert that ${x\qeu x'}\eqv {\bine x\reil U\bine{x'}}$ for 
all $x,\,x'\in \Om_{\log U}.$ 
(Consider any $x<x'$ in $\Om_{\log U}.$ 
If $d=x'-x\in U$ then $d<2^{a}$ for some $a\in\log U.$ 
As $\bine x\in\fse_{\log U},$ there is $b\in\log U\zt b>a,$ 
with $\bine x(b)=0.$ 
But easily $\bine x(z)=\bine{x'}(z)$ for any $z>b,$ hence, 
$\bine x\reil U\bine{x'}.$ 
The converse is obvious.) 

Yet for any $x\nin \Om_{\log U}$ there is 
$\tx\in\Om_{\log U}$ with $|x-\tx|\in U:$ 
put $\tx=x+2^{a+1},$ where $a$ is the largest number in 
$\log U$ with $\bine x(a)=0.$ 
For $x\in \Om_{\log U}$ put $\tx=x.$ 
The map $\vt(x)=\bine{\tx}$ is a Borel reduction of 
$\qeu$ to $\reil U.$

Finally, the map $f(\sg)=\sum_{z\in\adN}2^{2z}\sg(z)$ 
is a reduction of $\reil U$ to $\qeu.$ 
(The factor $2$ in $2z$ helps to avoid the   
trouble with values $\nin\Om_{\log U}.$) 
\epF

\brem
\label{rcof}
Choose $d\nin\log U.$ 
A slight modification of the same argument proves that 
$\qeu\boe {\Exp{d\bez U}}\ti\rav\adN.$ 
($d=[0,d),$ as usual.)
\erem


An obvious case when $\sis{a_n}{}\lek\sis{b_k}{}$ for 
increasing sequences is when 
$\frat{a_{n+1}}{a_n}\le\frat{b_{n+1}}{b_n}$ for all $n.$ 
The following result  
shows that this case essentially exhausts all cases of 
$\sis{a_n}{}\lek\sis{b_k}{}.$
Say that two increasing sequences 
$\sis{a_n}{}$ and $\sis{\al_n}{}$ are 
{\it cofinally equivalent\/} if $\cu{a_n}=\cu{\al_n}.$
Say that two decreasing sequences 
$\sis{a_n}{}$ and $\sis{\al_n}{}$ are 
{\it coinitially equivalent\/} if $\cv{a_n}=\cv{\al_n}.$ 

\bpr
\label{P:id}
Any two increasing sequences\/ $\sis{a_n}{}\zt\sis{b_k}{}$ 
are\/ \dd\lek compar\-able, in addition, 
if\/ $\sis{a_n}{}\lek\sis{b_k}{}$ then there are 
sequences\/ $\sis{\al_n}{}\zd\sis{\ba_k}{},$ cofinally 
equivalent to, resp., $\sis{a_n}{}\zd\sis{b_k}{},$ 
with\/ 
$\frat{\ba_{n+1}}{\ba_{n}}\ge\frat{\al_{n+1}}{\al_{n}}$ 
for all~$n.$

Similarly, any two decreasing sequences\/ 
$\sis{a_n}{}\zt\sis{b_k}{}$ 
are\/ \dd\lek comparable, in addition, 
if\/ $\sis{a_n}{}\lek\sis{b_k}{}$ then there 
are sequences\/ $\sis{\al_n}{}\zd\sis{\ba_k}{},$ coinitially 
equivalent to resp.\ $\sis{a_n}{}\zd\sis{b_k}{},$ 
such that\/ 
$\frat{\ba_{n}}{\ba_{n+1}}\ge\frat{\al_{n}}{\al_{n+1}}$ 
for all~$n.$
\epr
\bpF
We concentrate on the case of increasing sequences, the 
case of decreasing sequences is similar.
The conjunction of two symmetric forms of \eqref{nke+} is 
obviously contradictory, which implies the  
\dd\lek comparability assertion. 

Put $k_0=0$ and choose $n_0$ in accordance with
\eqref{ke+}, thus,
\addtocounter{equation}1
\envur{
\tag{$\mtho\arabic{equation}.0$}
\label{id2.0}
\kaz n'>n_0\;\sus k'>k_0\;
\big(\frakk{a_{n'}}{a_{n_0}}\le\frakk{b_{k'}}{b_{k_0}}\big).
}
If ({\sl Case 1\/}) we also have   
${\kaz k'>k_0\:\sus n'>n_0\:
\skl\frat{a_{n'}}{a_{n_0}}\ge\frat{b_{k'}}{b_{k_0}}\skp},$
then the sequences $\sis{a_{n_0+i}}{i\in\dN}$ and
$\sis{\frat{a_{n_0}+b_{k_0+i}}{b_{k_0}}}{i\in\dN}$
are cofinally equivalent, hence,
$\al_i=\frat{a_{n_0}+b_{k_0+i}}{b_{k_0}}$ and
$\ba_i=b_{k_0+i}$ prove the lemma.
Otherwise ({\it Case 2\/}) there is $k_1>k_0$ such that
$\frat{a_{n'}}{a_{n_0}}<\frat{b_{k_1}}{b_{k_0}}$ for all 
$n'>n_0.$
Choose $n_1>n_0$ so that, by \eqref{ke+}, 
\envur{
\tag{$\mtho\arabic{equation}.1$}
\label{id2.1}
\kaz n'>n_1\;\sus k'>k_1\;
\big(\frakk{a_{n'}}{a_{n_1}}\le\frakk{b_{k'}}{b_{k_1}}\big).
}
If we have now Case 1, \ie, symmetrically, 
$\kaz k'>k_1\;\sus n'>n_1\;
\skl\frat{a_{n'}}{a_{n_1}}\ge\frat{b_{k'}}{b_{k_1}}\skp,$
then, as above, the lemma holds immediately.
Thus, we can assume that there is $k_2>k_1$ with
$\frat{a_{n'}}{a_{n_1}}<\frat{b_{k_2}}{b_{k_1}}$ for all 
$n'>n_1.$
Choose $n_2>n_1$ following \eqref{ke+}.
{\it And so on\/}.

In the course of this construction, either the required 
result comes up just at some step, or we obtain increasing 
sequences $\sis{n_i}{}$ and $\sis{k_i}{}$ such that 
$\frakk{a_{n'}}{a_{n_i}}\le\frakk{b_{k_{i+1}}}{b_{k_i}}$ 
for all $n'>n_i$ and $i\in\dN.$
Let $\al_i=a_{n_i}\zt \ba_i=b_{k_i}$.
\epF

\seci{Countably cofinal monadic relations}
\label{ccof}

The goal of this section is to prove the part of 
\ref{T2} of Theorem~\ref T related to countably 
cofinal cuts and associated monadic \eqr s. 

Choose increasing sequences $\sis{a_n}{}\zt\sis{b_k}{}$ in 
$\adN$ with $U=\cud{a_n}$ and $V=\cud{b_k}.$ 
(Note that $\log U=\cu{a_n}$ and $\log V=\cu{b_k}.$) 
We are going to prove that $\qeu\qdr\qev$ iff 
$\qeu\bor\qev$ iff $U\lek V;$ 
the \dd\bor comparability of $\qeu\zd\qev$ then 
immediately follows from Proposition~\ref{P:id}.\vom 

{\sl Part 1\/}. 
Suppose that $\qeu\qdr\qev.$ 
Then $\reil U\qdr\reil V$ by Proposition~\ref{pcof}.
Let $\vt:\fse\to\fse$ be a CD reduction of $\reil U$ to 
$\reil V,$ thus, $\vpi\reil U \vpi'$ iff 
$\vt(\vpi)\reil V\vt(\vpi')$ for
all $\vpi,\,\vpi'\in\fse.$ 
The graph of $\vt$ has the form $\bigcup_{f\in F}C_f,$ 
where $F\sq\dn$ and $C_f=\bigcap_mC_{f\res m}$ for any 
$f\in\dn,$ sets $C_s\zt s\in\dln,$ are internal, and 
$C_t\sq C_s\sq\fse\ti\fse$ for $s\subset t,$ as in 
\eqref{cd} of Section \ref N.

Suppose, towards the contrary, that 
$\sis{a_n}{}\not\lek\sis{b_k}{},$ hence, we have 
\eqref{nke+}.

Suppose that $f\in F.$
Then $C_f$ is a subset of the graph of $\vt,$ hence,
by the choice of $\vt,$ for any $k\in\dN$ we have, for all 
$\vpi,\,\vpi',\,\psi,\,\psi'\in\fse,$ 
\dm
\kaz m\:\skl
{\otn\vpi{C_{f\res m}}{\psi}}\land 
{\otn{\vpi'}{C_{f\res m}}{\psi'}}
\skp\land
\resb\psi{b_k}=\resb{\psi'}{b_k}\limp
\sus n\:\skl
\resb\vpi{a_n}=\resb{\vpi'}{a_n}
\skp\,,
\dm
where $\resb\sg c=\sg\res(\adN\bez[0,c))$ for 
$\sg\in\fse$ and $c\in\adN.$
Then, by \Sat, 
\mul{
\label{eq1}
\kaz k\;\sus n\;\sus m\;
\kaz \vpi,\,\vpi',\,\psi,\,\psi'\in\fse:\\[0.4\dxii]
{\otn\vpi{C_{f\res m}}{\psi}}\land 
{\otn{\vpi'}{C_{f\res m}}{\psi'}}
\land \resb\psi{b_k}=\resb{\psi'}{b_k}\limp
\resb\vpi{a_n}=\resb{\vpi'}{a_n}\,.
}

\noi
A similar (symmetric) argument also yields the following:
\mul{
\label{eq2}
\kaz n\;\sus k\;\sus m\;
\vpi,\,\vpi',\,\psi,\,\psi'\in\fse:\\[0.4\dxii]
\otn\vpi{C_{f\res m}}{\psi}\land 
\otn{\vpi'}{C_{f\res m}}{\psi'}
\land \resb\vpi{a_n}=\resb{\vpi'}{a_n}
\limp \resb\psi{b_k}=\resb{\psi'}{b_k}\,.
}

To derive a contradiction to \eqref{nke+}, 
note first of all that  
$\cu{a_n}$ is a fast cut assuming \eqref{nke+}, thus,
we can suppose that $\frat{a_{n+1}}{a_n}$ is
infinitely large for all $n$ (Lemma~\ref{?fast}).
Now, let $k\in\dN$ witness \eqref{nke+}.
Let $n,\,m$ be numbers defined for this $k$ by \eqref{eq1}.
Choose $n'>n$ according to \eqref{nke+}: 
then $\frat{a_{n'}}{a_n}>\frat{b_{k'}}{b_k}$ for any $k'>k,$ 
hence, in fact, $\frat{a_{n'}}{a_n}>\ell+\frat{b_{k'}}{b_k}$ 
for any $m'>m$ and any $\ell\in\dN.$ 
Finally, choose $k'>k$ and $m'>m$ according to \eqref{eq2} 
but \wrt\ $n'.$
Put $C(f)=C_{f\res m'}.$
Then we have, for all 
$\ang{\vpi,\psi}\zt\ang{\vpi',\psi'}$ in $C(f):$
\pagebreak[0]%
\envur{
\label{><}
\left.
\bay{ccccc}
\resb\psi{b_k}=\resb{\psi'}{b_k}&\imp& 
\resb\vpi{a_n}=\resb{\vpi'}{a_n}&,&
\text{and}\\[0.7\dxii]
\resb\psi{b_{k'}}\ne\resb{\psi'}{b_{k'}}&\imp&
\resb\vpi{a_{n'}}\ne\resb{\vpi'}{a_{n'}}
\eay
\right\}\,.
}

We have $\fse=\dom\vt=\bigcup_{f\in F}X(f),$ 
hence, by \Sat, there is a finite set $F'\sq F$ such that 
still $\fse=\bigcup_{f\in F'}X(f).$ 
On the other hand, let us show that all sets $X(f)$ are 
too small for a finite union of them to cover $\fse.$ 
Call an internal set $X\sq\fse$ {\it small\/} iff 
\ben
\tenu{$\mtho(\fnsymbol{enumi})$}
\itla{smal}
there is a number $h\in\adN\bez\dN$ such that, for any 
internal map $\sg\in\fuk2{\adN\bez[0,h)}$ the set 
$X_\sg=\ans{\vpi\in X:\resb\vpi{h}=\tau}$ satisfies 
$2^{-h}\#X_\sg\simeq0.$ 
\een

\bpr
\label{PS}
$\fse$ is not a union of finitely many small 
internal sets.\qeD
\epr

It remains to show that any set $X(f)$ is small, with 
$h=a_{n'}$ in the notation above. 
(Note that $a_{n'}$ depends on $f,$ of course.) 
Take any $\ang{\vpi,\psi}\in C(f)$ and let 
$\sg=\resb\vpi{a_{n'}}\zt\tau=\resb\psi{b_{k'}}.$ 
By \eqref{><}, each $\ang{\vpi',\psi'}\in C(f)$ 
with $\resb{\vpi'}{a_{n'}}=\sg$ 
satisfies $\resb{\psi'}{b_{k'}}=\tau.$ 
Let us divide the domain 
$\Psi=\ans{\psi'\in\fse:\resb{\psi'}{b_{k'}}=\tau}$ 
onto subsets 
$\Psi_w=\ans{\psi'\in\Psi:\psi'\res{[b_k,b_{k'})}=w},$ 
where $w\in \fuk2{[b_k,b_{k'})}$ 
(\ie, $w$ is an internal map ${[b_k,b_{k'})}\to2$),
totally $2^{\frat{b_{k'}}{b_k}}$ of the sets 
$\Psi_w.$ 
For any such $\Psi_w,$ the set 
$\Phi_w=\ans{\vpi':
\sus \psi'\in \Psi_w\:\ang{\vpi',\psi'}\in C(f)}$ 
contains at most $2^{a_n}$ elements by the first 
implication in \eqref{><}, 
therefore, the whole set 
$X(f)_\sg=
\ans{\vpi'\in X(f):\resb{\vpi'}{a_{n'}}=\sg}$
contains at most $2^{a_n+b_{k'}-b_k}$ elements of the set 
$X(f),$ 
which is less than $2^{\frat{a_{n'}}\ell}$ for any 
$\ell\in\dN,$ hence, $X(f)$ is small, as required.\vom

{\sl Part 2\/}. 
Suppose that $\sis{a_n}{}\lek\sis{b_k}{},$ \ie, 
\eqref{ke+}, and derive $\reil U\bor\reil V.$ 
We can assume, by Proposition~\ref{P:id}, that 
$a_{n+1}-a_n\le b_{n+1}-b_n$ for all $n\in\dN.$
By Robinson's lemma, there is a number $N\in\adN\bez\dN$ and
internal extensions $\sis{a_\nu}{\nu\le N}$ and
$\sis{b_\nu}{\nu\le N}$ of sequences $\sis{a_n}{n\in\dN}$ and
$\sis{b_n}{n\in\dN},$ both being increasing hyperfinite
sequences satisfying $a_{\nu+1}-a_\nu\le b_{\nu+1}-b_\nu$
for all $\nu<N.$
Now we are ready to define a Borel reduction
$\vt$ of $\reil U$ to $\reil V$.

If $\vpi\in\fse$ then define $\vt(\vpi)=\psi\in\fse$ 
as follows:
\ben
\tenu{\arabic{enumi})}
\itla{xy1}
$\psi\res[0,b_0)$ is constant $0$ (not important);

\itla{xy2}
$\psi(b_\nu+h)=\vpi(a_\nu+h)$ whenever
$\nu<N$ and $h<a_{\nu+1}-a_\nu$;

\itla{xy3}
$\psi\res[b_\nu+a_{\nu+1}-a_\nu,b_{\nu+1})$
is constant $0$ for any $\nu<N$;

\itla{xy4}
$\psi(b_N+z)=\vpi(a_N+z)$ for all $z\in\adN$.
\een
Thus, to define $\psi,$ we move each piece
$\vpi\res[a_\nu,a_{\nu+1})$ of $\vpi$ so that
it begins with \dd{b_\nu}th position in $\psi,$ and
fill the rest of $[b_\nu,b_{\nu+1})$ by 0s;
in addition, $\psi\res[b_N,\iy)$ is a shift of
$\vpi\res[p_N,\iy).$ 
That $\vt$ is a Borel reduction of $\reil U$ to $\reil V$ 
is a matter of routine verification.

\seci{Countably coinitial monadic relations}
\label{ccoi}

That the double equivalence 
${\qeu\bor\qev}\eqv{\qeu\qdr\qev}\eqv{U\lek V}$ 
of \ref{T2} of Theorem~\ref T holds for any pair of 
countably coinitial cuts can be verified the same way as 
for countably cofinal cuts in Section~\ref{ccof} 
(with rather obvious amendments which account for the 
fact that now decreasing rather than increasing sequences 
$\sis{a_n}{}\zt\sis{b_k}{}$ are considered). 
We leave this to the reader, and concentrate, 
in this section, on \ref{T4} 
(the incomparability between countably cofinal and 
countably coinitial \er s), except for its Borel part. 

Suppose that $U=\cud{a_n}$ and $V=\cvd{b_k},$ where 
$\sis{a_n}{}$ and $\sis{b_k}{}$ are resp.\ (strictly) 
increasing and decreasing sequences of hyperintegers. 
Note that then $\log U=\cu{{a_n}}$ and 
$\log V=\cv{{b_k}}$.\vom

{\sl Part 1\/}.
Assuming that $\sis{a_n}{}$ is fast, prove that 
$\qeu\not\qdr\qev.$ 
We have a more general result: 
$\qeu\not\qdr{\qE}$ for any $\Pi^0_1$ equivalence 
relation $\qE$ on $\adN.$ 
It suffices (Proposition~\ref{pcof}) to show that 
$\reil U\not\qdr{\qE}.$
Suppose, towards the contrary, that $\vt:\fse\to\adN$ is 
a CD reduction of $\reil U$ to $\qE,$ so that 
${\vpi\reil U\vpi'}\eqv {\vt(\vpi)\qE \vt(\vpi')}$ 
for all $\vpi,\,\vpi'\in\fse.$
The graph of $\vt$ has the form $\bigcup_{f\in F}C_f,$ 
where $F\sq\dn$ and $C_f=\bigcap_mC_{f\res m}$ for any 
$f,$ all sets $C_s\zt s\in\dln,$ are internal, and 
$C_t\sq C_s\sq\fse\ti\adN$ whenever $s\subset t.$ 
Let ${\qE}=\bigcap_kE_k,$ where $E_k$ are internal sets 
and $E_{k+1}\sq E_k$ for all $k.$ 
As $\sis{a_n}{}$ is fast, we can assume that  
$\frat{a_{n+1}}{a_n}$ is infinitely large for any $n\in\dN$ 
(Lemma~\ref{?fast}).

By the choice of $\vt,$ for any $f\in F$ we have:
\mul{
\label{uv1}
\kaz \vpi,\,\vpi'\in\fse\;\kaz x,\,x'\in\adN:\quad 
\kaz m\:
(\otn\vpi{C_{f\res m}}x\land\otn{\vpi'}{C_{f\res m}}{x'})
\,\imp\,\\[0.4\dxii]
\skl
\sus n\:(\resb\vpi{a_n}=\resb{\vpi'}{a_n})\eqv 
\kaz k\:(\otn x{E_k}{x'})
\skp\,.
}

Applying \Sat\ here, with the implication $\mpi$ 
in the equivalence in the second line, we obtain 
numbers $m,\,n,\,k$ (which depend on $f$) such that  
\dm
{\otn\vpi{C_{f\res m}}x} \land
{\otn{\vpi'}{C_{f\res m}}{x'}}
\land{\otn x{E_{k}}x'}\limp 
\resb\vpi{a_{n}}=\resb{\vpi'}{a_{n}} 
\dm
for all $\vpi,\,\vpi'\in\fse$ and $x,\,x'\in\adN.$ 
Further, applying \Sat\ to \eqref{uv1} with the implication 
$\imp$ in the second line, with fixed numbers $k$ and $n+1,$ 
we find $m'(f)\ge m$ such that, for all 
$\vpi,\,\vpi'\in\fse$ and $x,\,x'\in\adN,$  
\dm
{\otn\vpi{C_{f\res m'(f)}}x}\land
{\otn{\vpi'}{C_{f\res m'(f)}}{x'}}\land 
{\resb\vpi{a_{n+1}}=\resb{\vpi'}{a_{n+1}}}\limp 
{\otn x{E_{k}}x'}.
\dm 
Let $X(f)=\dom C_{f\res m'(f)}.$  
It follows from the choice of $m'(f)$ that 
\dm
\kaz \vpi,\,\vpi'\in X(f):\quad 
{\resb\vpi{a_{n+1}}\ne\resb{\vpi'}{a_{n+1}}}
\lor {\resb\vpi{a_{n}}=\resb{\vpi'}{a_{n}}}\,,
\dm
therefore, $X(f)$ is small 
(see the definition before Proposition~\ref{PS}) 
because $\frat{a_{n(f)+1}}{a_{n(f)}}$ is infinitely large. 
This leads to a contradiction as in Section~\ref{ccof}.\vom

{\sl Part 2\/}. 
Prove that $\qev\not\qdr\qeu$ in any case. 
First of all, we can assume that $V$ is a slow 
countably coinitial cut, because if $V$ is such while 
$V'$ any countably coinitial cut then $\qev\qdr{\qei{V'}}$ 
by \ref{T2} of Theorem~\ref T. 
Thus, let $V=\cvd{d-k}=\bigcap_k[0,2^{d-k}),$ where 
$d\in\adN\bez\dN;$ then $\log V=\cv{b_k}=\bigcap_k[0,d-k).$ 
It suffices to prove that $\reil V\not\qdr\qeu$ 
(Proposition~\ref{pcof}). 
We show that, even more, $\reil V\not\qdr\qE$ for any 
$\Sg^0_1$ \eqr\ $\qE$ on $\adN.$
 
Suppose, on the contrary, that $\reil V\qdr\qE.$ 

Consider an auxiliary \eqr\ $\qR,$ defined on   
$\Xi=\fuk2d$ 
(all internal maps $d=[0,d)\to2$) 
as follows: $\sg\qR\tau$ iff 
$\sg\res{d\bez\log V}=\tau\res{d\bez\log V}.$~\footnote  
{\ 
Thus, $\qR$ is $\Exp{\ans{d-k:k\in\dN}},$ see Section~\ref{prel}, 
which is isomorphic to just $\Exp{\dN}$ on $\Xi$ via the 
bijection $\sis{i_z}{z<d}\mto\sis{i_{d-1-z}}{z<d}$ of $\Xi.$ 
In terms of this bijection, the partition of $\Xi$ into 
\dd\qR classes is equal to the partition onto  
into \dd dmonads  $M^d_g$ as in Section~\ref{sing}.}
For any $\sg\in \Xi$ let $\tilde\sg\in\fse$ be its extension 
by 0s. 
The map $\sg\to\tilde\sg$ is a reduction of $\qR$ to 
$\reil V,$ hence, in our assumptions, $\qR\qdr\qE.$ 
Let $\vt:\Xi\to\adN$ be a CD reduction of $\qE$ to $\reil U.$ 
Then $\vt=\bigcup_{f\in F}\bigcap_m C_{f\res m},$ where 
$F\sq\dn$ while $C_s\zt s\in\dln,$ are internal subsets of 
$\Xi\ti\adN$ with $C_s\sq C_t$ whenever $t\subset s.$ 
Finally, let ${\qE}=\bigcap_nE_n,$ where 
$E_n\sq\adN\ti\adN$ are internal sets and 
$E_n\sq E_{n+1}\zd\kaz n$.

For any $f\in F,$ we have, by the choice of $\vt,$ 
\mul{
\nonumber
\kaz \sg,\,\sg'\in \Xi\;\kaz x,\,x'\in\adN: 
\\[0.4\dxii]
\kaz m\:
({\otn\sg{C_{f\res m}}{x}}\land 
{\otn{\sg'}{C_{f\res m}}{x'}}
)
\land 
\kaz k\:(\resb\sg{d-k}=\resb{\sg'}{d-k})
\limp
\sus n\:(\otn x{E_n}{x'})\,,
}
where $\resb\sg{d-k}={\sg\res{[d-k,d)}}.$ 
Using \Sat, we obtain numbers 
$k=k(f)\zd n=n(f)\zd m=m(f)$ such that 
\envur{
\label{sgvp}
{\otn\sg{C_{f\res m}}{x}}
\land 
{\otn{\sg'}{C_{f\res m}}{x'}}
\land
{\resb\sg{d-k(f)}=\resb{\sg'}{d-k(f)}}
\limp
{\otn x{E_n}{x'}}\,.
}

We put $C(f)=C_{f\res m(f)}$ and $R(f)=\ran {C(f)}.$ 
It follows from \eqref{sgvp} that the set $R(f)$ can contain 
at most $2^{k(f)},$ a finite number, of pairwise 
\dd\qE inequivalent elements 
(because so is the number of all restrictions 
$\resb\sg{d-k(f)}\zd \sg\in \Xi$). 
On the other hand, since the graph of $\vt$ is covered by 
countably many sets of the form $C(f),$ the full image 
$\ran\vt=\ans{\vt(\sg):\sg\in \Xi}$ is covered by countably 
many sets of the form $R(f)$ 
(even if $F$ itself is uncountable), 
so that $\ran\vt$ contains only countably many pairwise 
\dd\qE inequivalent elements. 
Yet $\qR$ admits continuum-many 
pairwise \dd{\qR}inequivalent elements in $\Xi,$  
contradiction.

\seci{Remaining parts of the theorem on monadic \er s}
\label{tfin}

To check that $\rav\adN\bor\qeu$ for any additive \qd\  
cut $U,$ choose a number $c\nin U;$ then $x\mto xc$ is a 
Borel reduction of $\rav\adN$ to $\qeu,$  
in other words, $x=x'$ iff $xc\qeu x'c.$ 
This argument works for both countably cofinal and 
countably coinitial cuts $U.$ 

We continue with the following result, which proves the 
\dd\bor statement in \ref{T4} of Theorem~\ref T and ends 
the proof of \ref{T4} of Theorem~\ref T in general.

\ble
\label{NB}
If\/ $U$ is an additive countably cofinal cut and\/ $\qE$ a\/ 
$\Pi^0_1$ equivalence relation then\/ $\qeu\not\bor\qE$.  
\ele

It follows that\/ $\qeu\not\bor\qev$ provided\/ $V$ is any 
countably coinitial cut. 

\bpF
We know that $\dN\lek U$ (Remark~\ref{s.f}), hence, 
it can be assumed that $U=\dN.$ 
Let ${\qE}=\bigcap_nE_n,$ each $E_n\sq\adN$ internal and 
$E_{n+1}\sq E_n\zd\kaz n.$ 
Fix $c\in\adN\bez\dN$ and let $\vt:[0,c)\to\adN$ be a Borel 
reduction of ${\qen}\res[0,c)$ to $\qE.$ 
As any Borel (generally, any analytic) set, the graph of $\vt$ 
has the form $\bigcup_{f\in\nn}\bigcap_m C_{f\res m},$  
where $\nn$ is the set of all \dd\om sequences of natural 
numbers, all sets $C_u\sq[0,c)\ti\adN\zt u\in\nln,$ are internal, 
$\nln=$ all finite sequences of natural numbers, and 
$C_v\sq C_u$ whenever $u\subset v$ (see \cite{kkml}). 

Applying a simple measure-theoretic argument, we can find a 
sequence of numbers $\sis{j_m}{m\in\dN}$ in $\dN$ such that 
the set $X=\dom\vt'$ has Loeb measure $\ge \frac12,$ where 
$\vt'=\bigcup_{f\in F}\bigcap_m C_{f\res m}$ and 
$F=\ans{f\in\nn:\kaz m\:(f(m)\le j_m}.$ 
By Koenig's lemma, $\vt'=\bigcap_mC_m,$ where 
$C_m=\bigcup_uC_u,$ where the union is taken over all 
sequences $u$ of length $m$ such that $u(k)\le j_k$ for 
all $k<m,$ so that each $C_m$ is internal and (the graph of) 
$\vt'$ is a $\Pi^0_1$ set. 
Also, $\vt'={\vt\res X},$ where $X\sq[0,c)$ is a Borel set of 
Loeb measure $\ge \frac12$.

Since $\vt$ is a reduction (and $\vt'$ a partial one), we have
\mul{
\nonumber
\kaz x,\,x'\in X\;\kaz y,\,y'\in\adN:\\[0.3\dxii]
\kaz m\:
\skl
\otn x{C_m}y\land \otn {x'}{C_m}{y'})\limp
({\sus k\:(|x-x'|<k)}\eqv{\kaz n\:(\otn y{E_n}{y'}})
\skp\,.
}
Applying \Sat\ with $\mpi$ instead of $\eqv$ in the second line, 
we find numbers $m,\,n,\,k$ such that
\dm
\kaz x,\,x'\in X\;\kaz y,\,y'\in\adN:\quad
{\otn x{C_m}y}\land {\otn {x'}{C_m}{y'}}\land {\otn y{E_n}{y'}}
\limp {|x-x'|<k}\,.
\dm
Applying \Sat\ with $\imp$ instead of $\eqv,$ and fixed numbers 
$n$ and $4k,$ 
we find a number $m'\ge m$ such that
\dm
\kaz x,\,x'\in X\;\kaz y,\,y'\in\adN:\quad
{\otn x{C_{m'}}y}\land {\otn {x'}{C_{m'}}{y'}}\land {|x-x'|<4k}
\limp{\otn y{E_n}{y'}}\,.
\dm
It follows that ${|x-x'|<k}\lor{|x-x'|\ge4k}$ holds for all 
$x,\,x'\in X,$ which contradicts the assumption that $X$ has 
measure $\ge\frac12$.
\epF
 
\ref{T1} of Theorem~\ref T. 
Note that every CD-smooth \er\ is \dd\qdr reducible to 
$\qen$ because ${\rav\adN}\boe{\qen},$ see above. 
It follows, by \ref{T4} of Theorem~\ref T already proved, 
that $\qev$ is not CD-smooth (hence, not B-smooth), 
provided $V$ is an additive countably coinitial cut. 

Let $U=c\dN$ be a slow additive countably cofinal cut.
Note that $\qen$ has a \qd\ transversal $A$ by 
Theorem~\ref{sq}. 
Then $B=\ans{ac:a\in A}$ is obviously a CD transversal 
for $\qei{U},$ hence, $\qeu$ is CD-smooth 
(use the map sending any $x$ to the only element of 
$B$ equivalent to $x$). 
If a countably cofinal cut $U$ is fast then $U\not\lek\dN$ 
(say, by Lemma~\ref{?fast}), the 
non-CD-smoothness of $\qeu$ follows as above for countably 
coinitial cuts.

That $\qeu$ is not B-smooth for any additive 
countably cofinal cut $U$ follows from Lemma~\ref{NB}.
\vom

Finally, \ref{T5} of Theorem~\ref T.
It suffices, by \ref{T2}, to prove the following:

\ble
\label{den}
Suppose that, for any\/ $n,$ $\sis{a^n_k}{k\in\dN}$ is a 
fast increasing sequence. 
Then there are fast increasing sequences\/ $\sis{a_k}{}$ 
and\/ $\sis{b_k}{}$ such that\/ 
$\sis{a_k}{}\prec\sis{a^n_k}{k\in\dN}\prec\sis{b_k}{}$ 
for any\/ $n.$ 
The same for fast decreasing sequences.
\ele

Here $\sis{c_k}{}\prec\sis{d_k}{}$ means that 
$\sis{c_k}{}\lek\sis{d_k}{}$ but 
$\sis{d_k}{}\not\lek\sis{c_k}{}.$

\bpF
In the case of increasing sequences, 
we can assume that $d^n_k=\frat{a^n_{k+1}}{a^n_k}$ 
is infinitely large for all $n,\,k.$
By countable \Sat, there are numbers $a,\,b\in\adN\bez\dN$ 
such that $a<d^n_k<b$ for all $n,\,k.$
Put $a_k=k\sqrt a$ and $b_k=kb.$
\epF

\qeG{Theorem~\ref T}

\seci{An upper bound for countably cofinal relations}
\label{eo}

In classical descriptive set theory, the equivalence relation 
$\Eo,$ defined on $\dn$ so that $x\Eo y$ iff $x(n)=y(n)$ for 
all but finite $n,$ plays a distinguished role in the 
structure of Borel \er s, in particular, because it is the 
least, in the sense of Borel reducibility, non-smooth 
Borel equivalence relation. 
It would be a rather bold prediction to expect any analogous 
result in the ``nonstandard'' setting, yet a reasonable 
nonstandard version of $\Eo$ attracts some interest, giving a 
natural upper bound for countably cofinal monadic \er s. 

For $\xi,\,\eta\in\adS$ define: $\xi\FD\eta$ iff $\xi(x)=\eta(x)$ 
for all but finite $x\in\adN.$ 
($\FD$ from ``finite difference''.)

\ble
\label{Eo'}
If\/ $U\sq\adN$ is an additive countably cofinal cut then\/ 
$\qeu\bor\FD.$ 
If\/ $V\sq\adN$ is an additive countably coinitial cut then\/ 
$\qev\not\qdr\FD$.
\ele
\bpF
That $\qev\not\qdr\FD$ follows from the argument in 
Part~2 of Section~\ref{ccoi} because $\FD$ is obviously 
a $\Sg^0_1$ relation. 
As for the first statement, suppose that 
$U=\cud{a_n},$ where $\sis{a_n}{}$ is an increasing 
sequence in $\adN;$ accordingly, 
$\log U=\cu{a_n}=\bigcup_n[0,a_n).$ 
It suffices to prove that $\reil U\bor\FD.$ 

The sequence $\sis{a_n}{}$ admits an internal \dda extension 
$\sis{a_\nu}{\nu\le N},$ where $N\in\adN\bez\dN,$ still 
an increasing hypersequence of elements of $\adN.$  
Let, for any $\vpi\in\fse,$ $\vt(\vpi)$ be the 
(internal, hyperfinite) 
set of all restricted maps $\vpi\res[a_\nu,\iy)\zt \nu\le N,$ 
where $[a,\iy)=\adN\bez[0,a).$ 
By definition, $\vpi\reil U\psi$ iff the symmetric 
difference $\vt(\vpi)\sd\vt(\psi)$ is finite. 
Yet $\vt$ takes values in the set of all hyperfinite 
subsets of a certain internal hyper-countable set 
(because $\fse$ itself is hyper-countable) 
which can be identified with $\adN$.
\epF

\bcor
\label{eocor}
If\/ $U$ is as in the lemma then\/ $\qeu\bol\FD$. 
\ecor
\bpF
Use the lemma and \ref{T5} of Theorem~\ref T.
\epF

We don't know whether $\FD$ is an \poq{exact} upper bound for 
countably cofinal monadic \er s, but still the 
lower \dd\bor cone of $\FD$ contains many \er s not reducible 
to countably cofinal monadic ones, at least, all hyperfinite 
restrictions of $\FD$ are such. 
For any hyperfinite set $D\sq\adN$ let $\Fd D$ be the 
restriction of $\FD$ to the domain $\sui{(\fuk2D)},$ 
so that $\xi \Fd D\eta$ iff 
$\ans{d\in D:\xi(d)\ne\eta(d)}$ 
is finite. 
Easily ${\Fd D}\bor\FD,$   
moreover, ${\Fd D}\bol\FD$ because any possible CD 
reduction of $\FD$ to $\Fd D$ must be a bijection on 
any set $X\sq\fse$ of pairwise \dd\FD inequivalent elements, 
but we can take $X$ to be internal and hyper-infinite, 
which leads to contradiction because there is no CD 
injection from a hyper-infinite (internal) set in a 
hyperfinite set (say, by Lemma~\ref{carT}). 

\bte
\label{eo'}
If\/ $D$ is an infinite hyperfinite set and\/ $U$ an 
additive countably cofinal cut then\/ ${\Fd D}\not\qdr\qeu$.
\ete
\bpF
In the cource of the proof, it is more convenient to view 
$\Fd D$ as an equivalence on $\pv D$ defined so that 
$u \Fd D v$ iff $u\sd v$ is finite. 
Let, on the contrary, $\vt:\pv D\to\adN$ be a 
\qd\ reduction of $\Fd D$ to $\qeu.$ 
Assume that $D=[0,K)$ for some $K\in\adN\bez\dN.$ 
The graph of $\vt$ has the form 
$\bigcup_{f\in F}\bigcap_m P_{f\res m},$ where $F\sq\dn$ 
and $P_s\sq\pv D\ti\adN$ are as in \eqref{cd} of 
Section~\ref N. 
Let $X_s=\dom P_s$ and $X_f=\dom P_f,$ where 
$P_f=\bigcap_m P_{f\res m}$. 

Applying countable \Sat, we find a number 
$\nu\in\adN\bez\dN$ which is less than $K$ 
and moreover, $^\ast\!\rho(\nu)<K$ 
for any standard recursive function $\rho.$ 
Say that a set $Z\sq\pv D$ is {\it large\/} if there is 
an internal set $I\sq D$ such that $\#I=2\nu$ and 
$[I]^\nu\sq Z,$ where $[I]^\nu$ is the set of all internal 
subsets $Y\sq I$ with $\#Y=\nu.$
Then it is a consequence of the Ramsey theorem 
(in the nonstandard domain) 
that, for any $k\in\dN$ and any internal  
partition $\pv D=Z_1\cup...\cup Z_k$ 
at least one of the sets $Z_i$ is large.

We observe that there is $f\in F$ such that all sets 
$X_{f\res m}$ are large. 
(Otherwise let $X_{f\res m_f}$ be non-large for any 
$f\in F.$ 
Since $\dom\vt=\pv D,$ it follows from \Sat\ that 
$\pv D$ is a finite union of non-large sets of the form 
$X_{f\res m_f},$ contradiction with the above.) 
Then, by \Sat, $X_f$ itself is large, so that 
there is an internal set $I\sq X_f$ such that $\#I=2\nu$ 
and $[I]^\nu\sq X_f.$ 

Note that $P_f\sq\vt,$ hence, $P_f$ is a function, actually, 
$P_f={\vt\res {X_f}}.$ 
In addition, by \Sat, there is $n$ such that 
$\vpi=P_{f\res n}$ is already a function (internal). 
Then clearly $P_f=\vpi\res{X_f},$ therefore, 
${\vt\res{[I]^\nu}}={\vpi\res{[I]^\nu}},$ which implies that 
$\vt\res{[I]^\nu}$ is an internal map. 
Use this fact to derive a contradiction.

Let $I=\ans{a_1,...,a_{2\nu}}$ in the increasing order. 
For any $z=1,...,\nu,$ let 
$u_z=\ans{a_z,...,a_{z+\nu-1}}$ and 
$u_{\nu+z}=\ans{a_1,...,a_{z-1},a_{\nu+z},...,a_{2\nu}}$ 
(in particular, $u_{\nu+1}=\ans{a_{\nu+1},...,a_{2\nu}}$). 
Put $h_z=\vt(u_z).$ 
Easily the sets $u_z$ are internal and $\#{u_z}=\nu$ for all 
$z,$ moreover, $\#(u_z\sd u_{z+1})=2,$ hence, 
$u_z\Fd D {u_{z+1}}$ for each $z<2\nu,$ so that 
$|h_z-h_{z+1}|\in U$ because $\vt$ is a reduction, and, by 
the same reasons, $|h_{2\nu}-h_{1}|\in U.$ 
On the other hand, $\#(u_1\sd u_{\nu+1})=2\nu\nin\dN,$ hence, 
$|h_1-h_{\nu+1}|\nin U.$ 

To conclude, we have two hyperintegers 
$h_1$ and $h_{\nu+1},$ with $|h_1-h_{\nu+1}|\nin U,$ connected 
by two internal chains, $h_1,h_2,...,h_\nu,r_{\nu+1}$ and 
$h_{\nu+1},...,h_{2\nu},h_1,$ in which each link has length 
in $U.$ 
Obviously there is an index $z\zd1<z\le\nu,$ such that 
$|h_z-h_{\nu+z}|\in U.$ 
However by definition $\#(u_z\sd u_{\nu+z})=2\nu\nin\dN,$ 
hence, $|h_z-h_{\nu+z}|\nin U$ for any $z,$ contradiction.
\epF

Thus, we have the following two classes of \qd\ \eqr s 
strictly \dd\bor below $\FD:$ 
1)~\er s of the form $\qeu,$ where $U\sq\adN$ is an additive 
countably cofinal cut, 
2)~\er s of the form $\Fd {[0,c)},$ where $c\in\adN\bez\dN.$ 
It follows from our analysis that there is no \er\ in the 
first class \dd\qdr compatible with a \er\ in the second class. 
Is there anything below $\FD$ essentially different from 
these two classes~?

\seci{Final remarks and problems}
\label P

This final Section contains few scattered remarks and 
questions, 
mainly implied by analogies with \pol\ \dst.\vom

{\ubf Back to CD-cardinalities.}
Problem~\ref{?cuts} (Section~\ref{sing}) 
is, perhaps, the most interesting. 
Our analysis in the end of Section~\ref{sing} shows that, 
for $M^d_G$ to satisfy $M^d_G\qde c/\dN$ it is necessary and 
sufficient that $G\sq\dn$ is a set of Lebesgue measure $0.$ 
It is an interesting problem 
{\sl to find a reasonable necessary and sufficient 
condition for\/ $M^d_G$ to satisfy\/ 
$M^d_G\qde c\dN\qde [0,c).$\/} 
Can $M^d_G\qde [0,c)$ hold in the case when 
$\aleph_0<\card G<2^{\aleph_0}?$ 
Do these problems depend on the basic properties of 
the (standard) continuum in essential way~? 

How many \dd\boe classes of Borel subsets 
of $\adN$ do exist~?~\footnote
{\ This question can be addressed to \dd\qde classes of CD 
sets as well, but perhaps it is premature to search for 
an answer until Problem~\ref{?cuts} is solved.}
To answer such a question in the spirit of modern \dst, 
one has to define an \eqr, say, $\qE,$ on $\adN$ 
(in the \pol\ DST, on a Polish space), 
whose equivalence classes naturally represent \dd\boe 
classes of Borel subsets of $\adN,$ and classify 
$\qE$ in terms of best known, ``canonical'' \er s 
(see \cite{h,ndir}). 

It follows from Theorem~\ref{cardB} that Borel subsets of $\adN$ 
are represented, modulo $\boe,$ by sets of the following 
three classes: 
1) 
$\adN$ and cuts of the form $c=[0,c)\zt c\in \adN;$ 
2) 
additive countably cofinal cuts;
3) 
additive countably coinitial cuts. 

The first class naturally leads to ${\boe}\res\adN,$ \ie, 
the relation on $\adN$ defined so that $x\boe y$ iff 
there is a Borel bijection of 
$[0,x)$ onto $[0,y)$ iff $\frac xy\simeq 1.$  
Can it be characterized in terms of exponential equalities 
$\Exp D\zt D\sq\adN$? 
We conjecture that the relation ${\boe}\res\adN$ is 
\dd\boe equivalent to $\rav{\adN}\ti \Exp\dN.$ 

To approach the second class, fix $d\in\adN\bez\dN$ and 
let $D$ be the set of all increasing internal maps 
$\xi:d\to\adN$ 
satisfying $\xi(x+1)\ge x\xi(x)$ for all $x<d-1,$ so that 
any additive countably cofinal cut $U$ has the form 
$U=U(\xi)=\bigcap_{n\in\dN}\xi(n)$ for some 
(\poq{not} unique) $\xi\in D.$ 
Define $\xi\qE\eta$ iff $U(\xi)=U(\eta).$ 
This is a $\Pi^0_2$ equivalence; 
can it be described in terms of relations of the form 
$\rav X$ and $\Exp x$? 
Third class can be studied similarly, but with decreasing 
sequences and $\xi(x+1)\le \xi(x)/x$ for all $x,$ but 
does this lead to an \eqr\ \dd\boe equivalent to $\qE$?
\vom

{\ubf Exponential equalities.} 
Recall that $\Exp X$ is the \eqr\ of equality of internally 
extendable maps $X\to2,$ Section~\ref{prel}. 
This class of \er s contains, for instance, all monadic \er s 
(Proposition~\ref{pcof}, it suffices to take complements of CD cuts  
as sets $X$), 
hence, study of its properties in terms of $\qdr$ appears interesting 
and important. 
When ${\Exp X} \qdr {\Exp Y}$? 
The results for monadic \er s show that the answer has little to 
do with, for instance, the inclusion $X\sq Y.$ 
Our study of monadic \eqr s can be rather routinely generalized on 
\er s $\Exp X$ for sets $X\sq\adN$ of classes $\Sg^0_1$ and $\Pi^0_1$ 
(generalization of resp.\ countably coinitial and countably cofinal 
monadic \er s). 
For instance, it turns out that $\Exp X$ is not CD-smooth for 
any non-internal $\Sg^0_1$ set $X\sq\adN,$ as well as for any 
non-internal $\Pi^0_1$ set $X\sq\adN$ not of the form 
$H\bez C,$ where $H$ is internal and $C$ is countable. 
Is it true that $\Exp X$ is not CD-smooth for any set 
$X\sq\adN$ not in $\Pi^0_1$?
\vom

{\ubf A hyperfinite continuum-hypothesis.} 
Theorem~\ref{tCD} implies that, given $c\in\adN\bez\dN,$ 
there is no {\it regular\/} (see Section~\ref{sing}) 
CD-cardinalities strictly between those of $c/\dN$ and $c\dN$ 
(it is a question whether there are {\it singular\/} ones there). 
Are there any other similar pairs in the \dd\qdr structure~? 
A natural analogy with the continuum-hypothesis leads to the 
following question. 
Let $U$ be an additive CD cut in $\adN.$ 
(Or, generally, any CD subset of $\adN,$ but then the problem 
is most likely more difficult.)  
Does there exist any \qd\ \er\ $\qE$ with $\rav U\qdl\qE\qdl\Exp U$? 
Since $\rav U$ is the equality on $U$ while $\Exp U$ is the equality 
of internally extendable maps $U\to2,$ the double inequality can be 
seen to represent the fact that the CD-cardinality of the quotient 
space of $\qE$ is strictly between the CD-cardinality of $U$ and its 
natural ``power cardinality''. 
This question deserves a brief consideration. 

Let $d\in\adN\bez U,$ so that 
$\Exp U$ can be seen as the relation on $\fuk2d$ 
defined so that $\xi\Exp U \eta$ iff $\xi\res U=\eta\res U.$ 
That $\rav U\qdr\Exp U$ can be witnessed by the map 
$x\mto\xi_x,$ where $\xi_x\in\fuk2d$ is the characterictic 
function of the singleton $\ans x.$ 
If $U=H\bez C,$ where $H$ is internal while $C$ 
countable, then we can prove, using Lemma~\ref{carT}, that, 
paradoxically, $\rav U\qde\Exp U.$ 
Otherwise (see a remark above) $\Exp U$ is not CD-smooth, 
hence, $\rav U\qdl\Exp U$ strictly.  
Further, if there is $c\in U$ with $2^c\nin U$ then easily 
there are plenty of numbers $a<2^c\zt a\nin U$ with 
$\rav U\qdl\rav{a}\qdl\Exp U,$ thus, the 
``continuum-hypothesis'' fails. 

Now suppose that $U$ is exponentially closed, so that  
${c\in U}\imp{2^c\in U}.$ 
Then (Lemma~\ref{carT} applied) 
there is no CD set $X\sq\adN$ with 
$\rav U\qdl\rav{X}\qdl\Exp U,$ but is there any other 
\qd\ \er\ $\qE$ strictly \dd\qdr between  
$\rav U$ and $\Exp U$?\vom

{\ubf Another family of equivalence relations.} 
For any cut $\pu\ne U\sneq\adN,$ take $c\nin U$ and define, 
for (internal) $\xi,\,\eta\in\fuk2c,$ 
$\xi\qF_U\eta$ iff there are numbers $a\in U$ and 
$b\nin U \zd b\le c$ such that $\xi\res[a,b)=\eta\res[a,b).$ 
If $U$ is \qd\ then it belongs to $\Sg^0_1\cup\Pi^0_2,$ 
subsequently, $\qF_U$ can be transformed to $\Sg^0_2$ using 
\Sat. 
Anything about the \dd\qdr structure of this family~?\vom

{\ubf Smoothness and transversals.}
Our general method to establish smoothness was to find a 
suitable transversal.
Recall, in this context, that $\qen$ admits a \qd\ transversal by
Theorem~\ref{sq}, hence, is CD-smooth, but is not B-smooth
(Lemma~\ref{NB}), hence, does not admit a Borel transversal.   
However the existence of a transversal is not a necessary
condition for the smoothness.
Indeed, there exist Borel and B-smooth equivalence relations 
which do not admit even a \qd\ transversal~! 
An example can be easily extracted from the observation made 
in \cite[4.8]{kkml} that there is a $\Pi^0_2$ set in 
$\adN\ti\adN$ which does not admit a CD uniformization. 
\vom

{\ubf ``Fine structure'' of \eqr s.} 
Is the \er\ $\FD$ defined in Section~\ref{eo} in any sense 
\dd\qdr minimal over countably cofinal monadic \er s~?

Is there any result analogous to the Glimm -- Effros dichotomy 
(see \cite{hkl} or \cite{ndir}) 
of \pol\ \dst, in the same way as our Theorem~\ref{Ts} 
is analogous to the Silver -- Burgess dichotomy~? 
We conjecture that any Borel \eqr\ $\qE$ on $\adN$ is either 
B-smooth or satisfies ${\qen\res{[0,c)}}\bor\qE$ for some 
$c\in\adN\bez\dN$ or satisfies $\Exp\dN\bor\qE$. 

Theorem~\ref{sq} says that any countable CD \eqr\ is 
CD-smooth. 
What is the \dd\bor structure of countable {\it Borel\/} 
\er s~?\vom

{\ubf Ergodic theory.} 
Let $c\in\adN\bez\dN.$ 
The relation $\qen\res{[0,c)}$ on $[0,c)$ has certain 
similarities with the Vitali equivalence $x\qV y$ iff 
$x-y$ is rational on $\dR,$ for instance, Borel non-smoothness, 
the nonexistence of Borel transversals, perhaps, the 
\dd\bor minimality amongst all non-smooth \er s. 
However $\qen\res{[0,c)}$ lacks the following relevant 
property of $\qV:$ while every \dd\qV invariant Borel subset 
of $\dR$ has Lebesgue measure $0$ or its complement has 
measure $0,$ there exist plenty of \dd\qen invariant 
Borel subsets of $[0,c)$ 
having Loeb measure, for instance, $1/2:$ just take the 
\dd\qen saturation of $[0,\frac c2).$
(We consider the Loeb measure associated with the counting
measure $\mu(X)=\frac {\#X}c$ for internal subsets of
$[0,c).$)
Are there naturally defined ``nonstandard'' \er s which,
unlike $\qen,$ satisfy this property~?
Henson and Ross \cite[2.3]{hr} ask whether there exists a
bijection$f:[0,c)\onto[0,c)$ ergodic in the sense that for any
Loeb measurable set $X\sq[0,c)$ such that $X\sd f\ima X$ has
Loeb measure $0,$ the set $X$ itself has Loeb measure either
$0$ or $1;$ they prove that Borel bijections 
(\ie, with a Borel graph, as usual) are not ergodic. 
\vom

{\ubf Domain-independent version.} 
Define $\qE\bor' \qF$ if 
${\qE\ti{\rav\adN}}\bor{\qE\ti{\rav\adN}}.$ 
With this definition, we have, for instance, 
$\rav X\boe' \rav Y$ for any infinite hyperfinite $X,\,Y,$ 
and ${\qen\res a}\boe'{\qen\res b}\boe' \qen$ for any 
$a,\,b\in\adN\bez\dN,$ leading to structures less 
contaminated by the dependence on the size of the domain.

There is another possible way to the same goal.
Unlike the case of Polish spaces, it is not true in the
nonstandard domain that any Borel-measurable function
(\ie, here, it means that all preimages of internal sets
are Borel)
is Borel in the sense that its graph is Borel. 
It is known that, for rather good nonstandard universes,
for instance, those satisfying {\it the Isomorphism Property\/},
for any two infinite hyperfinite sets $X,\,Y$ there is a
bijection $f:X\onto Y$ such that the images and preimages of
internal sets are Borel.
(Such a bijection cannot be even \qd\ unless the fraction
$\frac{\#X}{\#Y}$ is neither infinitesimal nor infinitely
large.)
As mentioned in \cite{hr}, such a bijection induces an
isomorphism of the entire structure of Borel and \qd\ sets.

This naturally leads to the reducibility via Borel-measurable
maps.
Is $\qen$ Borel-measurable reducible to $\rav\adN$ in a
nonstandard universe satisfying the Isomorphism Property~?

\small

\let\section=\subsection


\begin{thebibliography}{99}

\bibitem{hkl}
L.~A.~Harrington, A.~S.~Kechris, A.~Louveau, 
A Glimm -- Effros dichotomy for Borel equivalence relations, 
{\it J.\ Amer.\ Math.\ Soc.\/} 1988, 310, pp.\ 293 -- 302.

\bibitem{he}
C.~W.~Henson,
Unbounded Loeb measures,
{\it Proc.\ Amer.\ Math.\ Soc.\/} 1979, 64, pp.\ 143 -- 160.

\bibitem{hr}
C.~W.~Henson and D.~Ross, 
Analytic mappings on hyperfinite sets,
{\it Proc.\ Amer.\ Math.\ Soc.\/} 1993, 118, pp.\ 587 -- 596.

\bibitem h
G.\ Hjorth,
{\it Classification and Orbit Equivalence Relations\/}
(Mathematical surveys and monographs, 75), AMS, 2000.

\bibitem{j}
R.~Jin,
Existence of some sparse sets of nonstandard natural 
numbers, 
{\it J. Symbolic Logic\/} 2001, 66(2), pp.\ 959 -- 973. 

\bibitem{dst}
A. S. Kechris.
{\it Classical Descriptive Set Theory\/}, 
Springer, 1995.

\bibitem{ndir}
A.~S. Kechris,
New directions in descriptive set theory,
{\it Bull. Symbolic Logic\/}, 1999, 5(2), pp.\ 161--174.

\bibitem{kkml}
H.~J.~Keisler, K.~Kunen, A.~Miller, and S.~Leth,
Descriptive set theory over hyperfinite sets,
{\it J. Symbolic Logic\/} 1989, 54, pp.\ 1167 -- 1180.

\bibitem{z}
B.~Zivaljevi\'c, 
Some results about Borel sets in descriptive set theory of 
hyperfinite sets, 
{\it J. Symbolic Logic\/} 1990, 55,  2, pp.\ 604--614.  


\end{thebibliography}
\end{document}